\documentclass[3p]{elsarticle}
\usepackage[usenames,dvipsnames]{color}
\usepackage{bm}
\usepackage{soul}
\usepackage[ampersand]{easylist}
\usepackage{amssymb,amsbsy,exscale,amsmath,amsfonts,amssymb,amscd}
\usepackage{graphicx}
\usepackage{mathtools}
\usepackage{stmaryrd}
\usepackage{bbm}
\graphicspath{{./images/}}
\usepackage{parskip}
\usepackage{longtable}
\usepackage{booktabs}
\usepackage{subcaption}
\usepackage{multirow}
\usepackage{array}
\usepackage[mathscr]{euscript}
\usepackage{algorithmicx}
\usepackage[ruled]{algorithm}
\usepackage{algpseudocode}
\usepackage{epstopdf}
\biboptions{sort&compress}
\usepackage{xspace}

\newcommand{\mbb}[1]{\mathbb{#1}}

\newcommand{\mcal}[1]{\mathcal{#1}}

\newcommand{\abs}[1]{\bigl\lvert{#1}\bigr\rvert}

\newcommand{\Linux}{\texttt{Linux}\xspace}

\newcommand{\refdom}{\Omega_{\bf r}}
\newcommand{\refdomB}[1]{\partial\Omega_{{\bf r},{#1}}}
\newcommand{\phdom}{\Omega}
\newcommand{\mcT}[1]{\mathcal{T}_h^{#1}}
\newcommand{\bx}{{\bf x}}

\newcommand{\bn}{{\bf n}}
\newcommand{\br}{{\bf r}}

\newcommand{\refN}{{\bf n}_{\bf r}}
\newcommand{\refX}{{\bf x}_{\bf r}}
\newcommand{\wn}{\kappa^2}
\newcommand{\nX}{\nabla_{\bf r}}

\newcommand{\bz}{{\bm z}}

\providecommand{\abs}[1]{\lvert#1\rvert}
\providecommand{\norm}[1]{\lVert#1\rVert}

\newcommand{\lp}{\left(}
\newcommand{\rp}{\right)}

\newcommand{\lb}{\left[}
\newcommand{\rb}{\right]}

\newcommand{\beql}{\begin{equation}}
\newcommand{\eeql}{\end{equation}}
\newcommand{\beq}{\begin{equation*}}
\newcommand{\eeq}{\end{equation*}}
\newcommand{\bali}{\begin{align}}
\newcommand{\eali}{\end{align}}
\newcommand{\balit}{\begin{alignat}}
\newcommand{\ealit}{\end{alignat}}
\newcommand{\bse}{\begin{subequations}}
\newcommand{\ese}{\end{subequations}}
\newcommand{\bfig}{\begin{figure}[!ht]}
\newcommand{\efig}{\end{figure}}

\newcommand{\eMat}{\end{matrix} \right]}

\providecommand{\norm}[1]{\lVert#1\rVert}

\providecommand{\apo}{\emph{a posteriori }}

\newcommand{\Ff}{\mathcal{F}_{\mathsf{f}}}

\newcommand{\MATLAB}{\textsc{Matlab}\xspace}

\newcommand\RT[1]{\textcolor{black}{#1}}
\def\etal{{\it et al. }}

\journal{Journal of Computational Physics}

\begin{document}

\begin{frontmatter}

\title{A multiscale continuous Galerkin method for stochastic simulation and robust design of photonic crystals}

\author[mit]{F.~Vidal-Codina\corref{cor1}\fnref{fn1}}
\ead{fvidal@mit.edu}
\author[mit]{J.~Sa\`{a}-Seoane\fnref{fn1}}
\ead{joel.saa@gmail.com}
\author[mit]{N.-C.~Nguyen\fnref{fn1}}
\ead{cuongng@mit.edu}
\author[mit]{J.~Peraire\fnref{fn1}}
\ead{peraire@mit.edu}

\cortext[cor1]{Corresponding author}
 \address[mit]{Department of Aeronautics and Astronautics, Massachusetts Institute of Technology, Cambridge, MA 02139, USA}
\fntext[fn1]{This work was supported by AFOSR Grant FA9550-15-1-0276.}


 \begin{abstract}
We present a multiscale continuous Galerkin (MSCG) method for the fast and accurate stochastic simulation and optimization of time-harmonic wave propagation through photonic crystals. The MSCG method exploits repeated patterns in the geometry to drastically decrease computational cost and incorporates the following ingredients: (1) a reference domain formulation that allows us to treat geometric variability resulting from manufacturing uncertainties; (2) a reduced basis approximation to solve the parametrized local subproblems; (3)  a gradient computation of the objective function; and (4) a model and variance reduction technique that enables the accelerated computation of statistical outputs by exploiting the statistical correlation between the MSCG solution and the reduced basis approximation. The proposed method is thus well suited for both deterministic and stochastic simulations, as well as robust design of photonic crystals. We provide convergence and cost analysis of the MSCG method, as well as a simulation results for a waveguide T-splitter and a Z-bend to illustrate its advantages for stochastic simulation and robust design.
 \end{abstract}
\begin{keyword}
Photonics \sep multiscale methods \sep continuous Galerkin \sep model reduction \sep variance reduction  \sep robust design
  \end{keyword}
  
\end{frontmatter}

\section{Introduction}

A field that has attracted significant interest in recent years is the study of light propagation through photonic crystals, which are heterogeneous materials engineered to exhibit properties that cannot be found in homogeneous materials. Photonic crystals are assembled by combining conventional materials in lattice structures usually at the microscopic level. Electromagnetic wave propagation through photonic crystals is characterized by band gaps, that is,  ranges of frequencies for which light waves are not allowed to propagate through the periodic optical nanostructure \cite{saleh1991fundamentals}. Introducing defects in these periodic structures allows us to create waveguides, since waves traveling at frequencies within the band gap are exponentially attenuated within the crystal, and thus can only propagate along the defect.  Waveguides exploiting this mechanism can be much more efficient than traditional waveguides based on total internal reflection (TIR). Photonic crystals have applications in fibers \cite{knight1996all,russell2003photonic}, waveguides \cite{johnson1999guided,joannopoulos2011photonic} and superlenses \cite{pendry2000negative,luo2002all}.


The simulation of wave propagation phenomena in heterogeneous media is an active research field and several techniques have been proposed. One of the most widely used approaches is the finite-difference time-domain (FDTD) method  \cite{taflove2005computational,kunz1993finite,sullivan2013electromagnetic,oskooi2010meep}. The FDTD method is simple and efficient, albeit not well suited to treat complex geometries, irregular domains and multiple length scales.  Another family of approaches are the finite-volume time-domain (FVTD) methods \cite{balsar2017computational,balsar2018computational}, which enable mesh adaptation and refinement and exhibit low dissipation and dispersion due to the use of high-order Godunov schemes. The finite element (FE) method both in the time and frequency domain  \cite{sadiku2000numerical,jin2014finite,yasumoto2005electromagnetic} is also a popular alternative for solving wave propagation problems thanks to its ability to handle complex geometries and inhomogeneous materials, as well as allow $h/p$ adaptivity. Furthermore, material interface conditions and boundary conditions can be implemented in a natural manner. Although low order methods are often used due to their simplicity,  high order methods \cite{graglia1997higher,ainsworth2003hierarchic} are more accurate and efficient if high accuracy is required. All the above methods have been used to study the propagation of waves in photonic crystals \cite{sigalas1993photonic,lidorikis2000gap,frei2004finite,rodriguez2005disorder,ismagilov2015second,hammerschmidt2016reconstruction}.

The ability to accurately simulate wave propagation in photonic crystals presents  unique challenges. Problems of interest typically  involve complex geometries and a mismatch in critical length scales, which can be of several orders of magnitude. Resolving the small scales using  uniform discretizations  requires a prohibitive large number of grid points.  An additional difficulty stems from the fact that mathematical models (e.g., based on Helmholtz and Maxwell's equations) may not capture the real physical phenomena accurately enough due to simplifications and uncertainty in the model data, such as geometry, material properties, and boundary conditions.

Several methods have been developed to deal with multiple scales in composite materials \cite{kanoute2009multiscale}. Homogenization methods \cite{bensoussan2011asymptotic} allow the treatment of multiscale features by solving a coarse-scale model which has been modified to account for the small scales. These approaches can be quite successful at predicting the global macroscopic behavior but lack the detailed description of the physics at the smaller scales. The multiscale FE method \cite{hou1997multiscale,hou1999convergence} and the mixed multiscale FE method \cite{chen2003mixed} are also alternatives that have been successfully applied to multiscale elliptic problems. {The main idea is to construct multiscale FE basis functions on a coarse-scale that capture the fine-scale features of the solution. These multiscale functions are then coupled globally into a linear system, whose unknowns are the solution values on the coarse grid. The fine-scale solution can be recovered combining the multiscale functions with the coarse-grid nodal values.} The main drawback of these methods is their strong dependence on the boundary conditions of the subproblems and, for problems where strong heterogeneities are present at the interfaces, one needs to develop adaptive boundary conditions to avoid small scale resonances \cite{hou1997multiscale}.  

In our approach, we advocate for not modifying  the original multiscale problem, but instead obtaining computational efficiency via domain decomposition.  A successful domain decomposition approach is the mortar element method \cite{bernardi1993domain, bernardi1994new}. This method allows for an independent discretization of each subdomain and enforces continuity of the solution across subdomains weakly.  Since the  meshes across the subdomain boundaries do not need to match, this approach allows great flexibility in the definition of the subdomains, and becomes particularly attractive for problems requiring mesh adaptivity and complex geometries. Other methods include the multiscale DG method \cite{aarnes2005multiscale}, which blends the multiscale FE method and imposes weak continuity at the subdomain interfaces, the hybridized multiscale DG method \cite{nguyen2013hybridized}, the geometric multiscale FEM  \cite{casadei2013geometric}, and the method of polarized traces \cite{zepeda2016method} used for high-frequency problems. For an extensive review of multiscale FE methods, we refer the reader to \cite{efendiev2009multiscale}.


In this paper, we present the multiscale continuous Galerkin (MSCG) method to simulate wave propagation problems on structured materials. This method is an extension of the hybridized continuous Galerkin method (CG) introduced in \cite{cockburn2007locally} and the hybridized multiscale DG method \cite{nguyen2013hybridized}. The multiscale CG method possesses considerable advantages over other simulation methods. First, the multiscale discretization results in a significantly smaller global linear system than that of the original problem, owing to static condensation of all the degrees of freedom corresponding to the domain interiors. Second, it exploits repetitive patterns in the structure to rapidly construct the global linear system by only solving a small number of local subproblems.  As a result, the method provides fast and accurate simulations of wave propagation in structured media beyond the capabilities of current numerical methods.

Furthermore, we also consider the simulation and optimization of wave propagation  problems in the presence of fabrication uncertainties. For complicated structures the interaction of electromagnetic waves with heterogeneous materials can be very sensitive to geometry, thus small perturbations in the microscopic structure may significantly degrade the performance of the device. To that end, we  augment the MSCG method with the following ingredients: (1) a reference domain formulation \cite{persson2009discontinuous} that allows us to treat geometric variability; (2) a reduced basis approximation  \cite{barrault2004empirical,prud2002reliable,rozza2008reduced} at the subdomain level \cite{phuong2013static,huynh2013static} that drastically accelerates the solution of the parametrized local subproblems; (3)  an adjoint technique for computing gradients of an output functional; and (4) a model and variance reduction technique \cite{vidal2015model,vidal2016empirical} that enables the accelerated computation of statistical outputs by exploiting the statistical correlation between the MSCG solution and the reduced basis approximation.



This article is organized as follows. In Section \ref{sec:MSCG}, we introduce the wave propagation problem and present the MSCG method. In Section \ref{sec:sMSCG}, we extend our approach to stochastic simulation and robust design problems. In Section \ref{sec:Res}, we present numerical results to demonstrate the performance of the proposed method. Finally, we summarize the main conclusions of this article in Section \ref{sec:conclusions}.

\section{The multiscale continuous Galerkin method}\label{sec:MSCG}
\subsection{Problem statement}
The propagation of the time-harmonic electric and magnetic fields ${\bf E},\,{\bf H}$ is governed by the frequency-domain Maxwell's equations. However, for two-dimensional problems the linearity of Maxwell's equations implies the transverse magnetic (TM) and transverse electric (TE) polarizations decouple, hence it is only necessary to solve the scalar Helmholtz equation in a domain $\Omega \in \mathbb{R}^2$, with Lipschitz boundary \RT{$\partial \Omega = \partial \Omega_D \cup \partial \Omega_N$ where Dirichlet and Neumann conditions are both prescribed}, namely:
\bse\label{eq:helmholtz}
\begin{alignat}{2}
 -\nabla \cdot \lp \rho(\bx)\nabla u \rp - \wn(\bx) u &= f,\qquad &\bx & \in \Omega\;, \label{eq:helmeq}\\
\rho \nabla u \cdot \bn  &= h ,\qquad&  \bx & \in \partial\Omega_{N} \;, \label{eq:helmN}\\
u\  &= u_D ,\qquad&  \bx & \in \partial\Omega_D \,. \label{eq:helmD}
\end{alignat}
\ese
\RT{The auxiliary parameters $\rho$ and $\kappa$ enable the compact definition of the above non-dimensional equation for both polarizations. In TM or E-polarization,  the magnetic field is confined to the plane of propagation (which we assume the $x-y$ plane), that is ${\bf H} = (H_x,H_y,0)$, and  the electric field is perpendicular to this plane ${\bf E} = (0,0,E_z)$, hence $\rho = 1$, $\wn = (\omega/c)^2\varepsilon(\bx)$ and $u = E_z$. Conversely, in TE or H-polarization, the electric field is confined to the plane and the magnetic field is perpendicular to it, that is ${\bf E} = (E_x,E_y,0)$ and ${\bf H} = (0,0,H_z)$, hence $\rho = \varepsilon(\bx)^{-1}$, $\wn = (\omega/c)^2$ and $u = H_z$. The dimensionless relative permittivity field $\varepsilon(\bx)$ is a piecewise-constant function that prescribes distinct permittivity values to different spatial regions within the photonic crystal, thus enabling the simulation of heterogeneous materials. The angular frequency of the wave $\omega$ is non-dimensionalized with the speed of light $c$ and the periodicity of the photonic crystal $a$ as $\overline{\omega} = \omega a/(2\pi c)$.}

The Helmholtz equation above assumes a finite computational domain. Nonetheless, wave propagation problems often occur in unbounded domains. In order to simulate unboundedness while restricting our simulation to a finite domain, we resort to Perfectly Matched Layers (PMLs) \cite{berenger1994perfectly} surrounding the photonic crystal. PMLs can be encoded in the auxiliary parameters with the transformations
\begin{equation}\label{eq:pml}
\rho \mapsto \left[ {\begin{array}{cc}
   \dfrac{s_y}{s_x} & 0 \\
   0 & \dfrac{s_x}{s_y} \\ 
  \end{array} } \right] \rho\;,\qquad \kappa^2 \mapsto s_xs_y\kappa^2\;,
\end{equation}
where $s_x$ is a complex-valued frequency-dependent damping function that prescribes the attenuation of the $x$-propagating waves inside the PML, and analogously for $s_y$. We refer to \cite{johnson2008notes} for a detailed description of PMLs, their numerical implementation and definition of damping functions.

\RT{\subsection{Structured heterogeneous materials}}

\begin{figure}[h!]
 \centering
 \includegraphics[scale = .65]{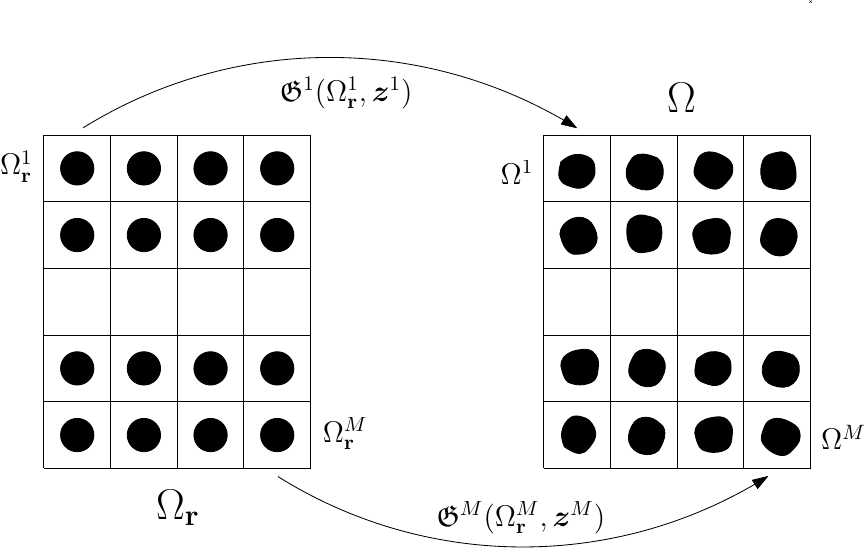}
 \caption{Reference domain (left), physical domain (right) and subdomain mapping.}\label{fig:motivation}
\end{figure}

As mentioned earlier, our goal in this paper is to develop a method to efficiently simulate wave propagation through heterogeneous structured materials. Heterogeneity is encoded in $\varepsilon(\bx)$ by assigning different permittivity values to distinct spatial regions forming the photonic crystal. As a motivating example, we consider the waveguide shown in Fig. \ref{fig:motivation} (left) consisting of a regular arrangement of dielectric rods $(\varepsilon(\bx) > 1)$ in air $(\varepsilon(\bx) = 1)$, and a row defect without rods. Our work relies on decomposing the domain into smaller subdomains that belong to distinct classes, hence within each class, the geometry of the subdomains is the same.  In this case, only two different classes of subdomains are needed: the subdomains with the rod and subdomains without the rod. Unfortunately, in real photonic applications the rods are not identical, mainly due to manufacturing variability that arise as a consequence of the extreme-scale fabrication techniques as sketched in Fig. \ref{fig:motivation} (right). Thus, in order to achieve truly predictive simulations one must take into account these variations. In this paper we consider only circular dielectric rods, but the extension to other types/shapes of rods is straightforward.

\RT{\subsection{Reference domain formulation}}

We assume that the geometry of the physical domain $\Omega$ is parametrized by $\bm z$, defined on a compact set, and that we want to solve \eqref{eq:helmholtz} for many different realizations of $\bm z$. In this scenario, it is much more convenient to map the physical domain $\Omega$ onto a fixed reference domain $\refdom$. For the waveguide example in Fig. \ref{fig:motivation}, the reference domain $\refdom$ is the structure with perfectly circular rods in Fig. \ref{fig:motivation} (left), whereas the physical domain $\Omega$ is the real waveguide with irregular rods in Fig. \ref{fig:motivation} (right). Furthermore, we can leverage the structure of the photonic crystal partitioning both the physical and the reference domain into $M$ non-overlapping subdomains, that is $\overline{\Omega} = \bigcup_{m=1}^M \overline{\Omega}^m$, $\overline{\Omega}_\br = \bigcup_{m=1}^M \overline{\Omega}_\br^m$. These partitions enable us to define $M$ independent diffeomorphisms $\lbrace\mathfrak{G}^m \rbrace_{m=1}^M$ parametrized by $\lbrace\bz^m \rbrace_{m=1}^M$  to prescribe geometry transformations at the subdomain level $\Omega^m = \mathfrak{G}^m(\refdom^m,\bz^m)$, while leaving the subdomain interfaces fixed. The overall mapping is then given by $\mathfrak{G} \coloneqq\mathfrak{G}^1\times\ldots\times\mathfrak{G}^M$ and $\bm z \coloneqq \bm z^1\times\ldots\times \bz^M $. 

Solving \eqref{eq:helmholtz} for a given $\bz$ first requires either re-discretizing the transformed geometry or  deforming the original discretization, which can be complex and cumbersome. Alternatively, we follow Persson \etal \cite{persson2009discontinuous} and use the diffeomorphism $\mathfrak{G}$ that maps $\refX \in \refdom$ to $\bx\in\Omega$ to transform the governing equations \eqref{eq:helmholtz} into a new set of equations defined on the reference domain $\refdom$ that remains unchanged. To derive these equations,  we first integrate on a volume and use the divergence theorem. The mapping $\mathfrak{G}$ has deformation gradient $\mcal{G} = \nX \mathfrak{G}$ and Jacobian $g = \det \mcal{G}$, and for an arbitrary control volume $v\subset\Omega$ (with Lipschitz boundary) and its associated $v_{\bf r}\subset\refdom$, the transformed equation on the reference domain reads
\begin{equation*}
\begin{split}
  0 = \int_v \lb \nabla\cdot  \rho\nabla u  +\wn u -f \rb dv & = \int_s   \rho\nabla u \cdot \bn ds + \int_v (\wn u -f) dv\\ &= \int_{s_{\bf r}} g\mcal{G}^{-1}  \rho\nabla u \cdot \refN d{s_{\bf r}}
  + \int_{v_{\bf r}} (\wn u-f) gdv_{\bf r} \\ &=  \int_{v_{\bf r}} \lb\nX \cdot g\mcal{G}^{-1}\rho\nabla u +\wn g u - fg  \rb dv_{\bf r} \\ &= \int_{v_{\bf r}} \lb\nX \cdot \rho g \mcal{G}^{-1}\mcal{G}^{-T} \nX u  +\wn g u - fg  \rb dv_{\bf r} \;,
\end{split}
\end{equation*}
where we have invoked the relationships $\nabla = \mcal{G}^{-T}\nX$, $dv = gdv_{\bf r}$ and $\bn ds = g \mcal{G}^{-1} \refN d{s_{\bf r}}$ with outward-pointing normals $\bn,\, \refN$. Transforming the boundary conditions, we arrive at
\bse\label{eq:helmholtz_ref}
\begin{alignat}{2}
 -\nX \cdot \lp \rho {\bm G} \nX u \rp - \wn u\,g &= f\,g,\qquad &\bx & \in \refdom\;, \label{eq:helmeqref}\\
\rho \nX u \cdot \refN  &= h\,g_s ,\qquad&  \bx & \in \partial \Omega_{\br,N} \;, \\
u\  &= u_D ,\qquad&  \bx & \in \partial \Omega_{\br,D} \;, 
\end{alignat}
\ese
where $\bm G = g \mcal{G}^{-1}\mcal{G}^{-T}$ and $g_s$ is the restriction of the Jacobian $g$ to the face.  Thus, for a given $\bz$ instead of solving \eqref{eq:helmholtz} on the deformed domain $\Omega$, we may alternatively solve \eqref{eq:helmholtz_ref} on the original domain $\refdom$.



\subsection{MSCG method}

\begin{figure}[h!]
\centering
\subcaptionbox{Domain partition into 5 subdomains (3 classes), with faces of subdomain 5 highlighted in dashed red. Subdomain interfaces are discretized with a single fourth order element except $\mathsf{f}^1$, where two elements $\mathsf{f}^1_1,\,\mathsf{f}^1_2$ of order four are used. \label{fig:mscg}}
[10cm]{\includegraphics[scale = 1]{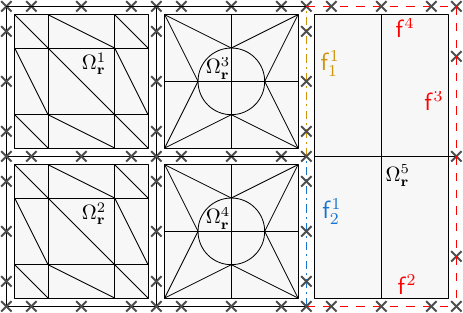}}
\hfill\subcaptionbox{Basis of Lagrange polynomial at subdomain boundary, using one element of order five per interface. \label{fig:globfun}}
[5cm]{\includegraphics[scale = 0.23]{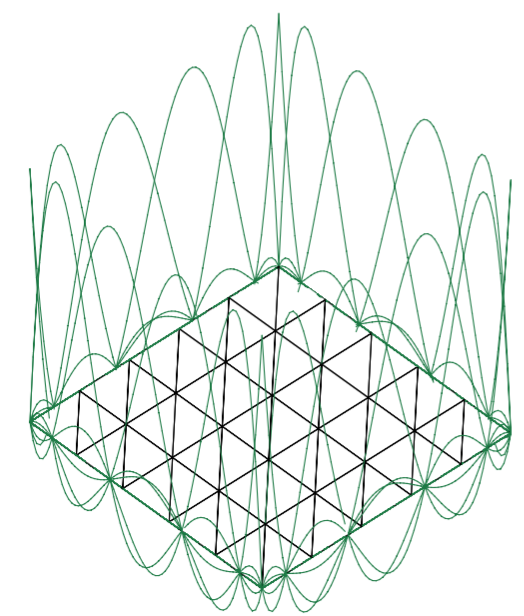}}
\caption{Multiscale schematics for local problems and Lagrange polynomial at interfaces.}
\end{figure}

Let us introduce the MSCG method to solve the transformed governing equations \eqref{eq:helmholtz_ref} on the fixed reference domain $\refdom$. We first define a discretization $\mcT{m}$ on each subdomain $\refdom^m$, consisting of disjoint regular elements $T$ that partition the subdomain. The discretization of the entire reference domain is thus $\mcT{} = \bigcup_{m=1}^M \mcT{m}$. Subdomain discretizations are disjoint between subdomains, thus one can seek for instance a triangulation for one subdomain and a quad discretization for its neighboring subdomain. The division of a sample reference domain into five subdomains with its corresponding discretizations is illustrated in Fig. \ref{fig:mscg}. 



Let  $\mathsf{f}^\ell$ be a subdomain face given by either $\partial \refdom^m \cap \partial \refdom^n$ ($n \neq m$) or $\partial \refdom^m \cap \partial \refdom$. The collection of subdomain faces of $\mcT{}$ is $\lbrace \mathsf{f}^\ell \rbrace_{\ell= 1}^L$. Furthermore, each subdomain face $\mathsf{f}^\ell$ may be subdivided into $N^\ell$ elements $\mathsf{f}^\ell_i,\, 1\le i \le N^\ell$. Thus we define the set of global face elements as $\Ff = \lbrace \mathsf{f}_i^{\ell},\, 1 \le i \le N^\ell, 1 \le \ell \le L \rbrace$, also referred to as the mesh skeleton. We note that, not only the discretizations of the faces and the boundaries of the adjacent subdomains are not required to coincide, but there is also flexibility in the way we define the subdomain faces. For instance, in Fig. \ref{fig:mscg},  the subdomain faces coincide with the global face elements except for the interface between $\refdom^5$ and both $\refdom^3,\,\refdom^4$, wherein the subdomain face $\mathsf{f}^1$ is subdivided into two global face elements $\mathsf{f}^1_1$ (gold) and $\mathsf{f}^1_2$ (blue). 

We can now introduce our approximation spaces as
\begin{align*}
W_h &= \{w\in L^2(\refdom) : w \in \mcal{C}^0(\refdom^m),\,w|_T \in \mcal{P}^{p^m}(T), \,\forall T\in \mathcal{T}^m_h,\, 1\le m \le M \},\\
M_h &= \{\mu\in L^2(\partial \refdom^m) : \mu|_{\partial \refdom^m}  = w|_{\partial \refdom^m},\, \mbox{ for } w \in W_h \},\\
V_h &= \{v \in \mcal{C}^0(\Ff) : v|_{\mathsf{f}} \in \mcal{P}^{p^\mathsf{f}}(\mathsf{f}),\, \forall \mathsf{f} \in \Ff \} ,
\end{align*}
where $\mcal{P}^{p^m}(T)$ is the space of polynomials of degree at most $p^m$ on $T\in \mcT{m}$ and $\mcal{P}^{p^\mathsf{f}}(\mathsf{f})$ is the space of polynomials of degree at most $p^\mathsf{f}$ on $\mathsf{f} \in \Ff$. Note that we allow for polynomial spaces of different degrees for both different subdomains and face elements on the subdomain's boundaries. To impose Dirichlet boundary conditions, we introduce  the set $V_h(u_D) = \{v \in V_h : v = \texttt{P}^{V_h}(u_D),\; \textnormal{on } \partial \Omega_{\br,D}\}$, where the operator $\texttt{P}^{V_h}$ represents the $L^2$ projection onto the space $V_h$ on the  boundary $\partial \Omega_{\br,D}$. \RT{The volume inner products for these finite element spaces are defined as
\begin{equation*}
 (\eta,\zeta)_{\refdom} :=  \sum_{m = 1}^M(\eta,\zeta)_{\refdom^m} =  \sum_{m = 1}^M\sum_{T\in\mcT{m}}(\eta,\zeta)_{T} = \sum_{m = 1}^M\sum_{T\in\mcT{m}}\int_T \eta\zeta\,,
\end{equation*}
and the surface inner products as
\begin{equation*}
 \langle \eta,\zeta\rangle_{\partial\refdom^m} :=  \sum_{T\in\mcT{m}} \langle \eta,\zeta\rangle_{\partial T} =  \sum_{T\in\mcT{m}}\int_{\partial T} \eta\zeta\;.
\end{equation*}}


We then introduce the auxiliary variable $q_h$ which approximates the normal component of the flux $q = \rho {\bm G} \nX u \cdot \refN$, and seek an approximation $(u_h,\lambda_h,q_h) \, \in W_h \times V_h(u_D) \times M_h$ such that
\begin{align*}
(\rho {\bm G} \nX u_h,\nX w)_{\refdom} -(g\,\wn u_h,w)_{\refdom} - \sum_{m=1}^M \langle q_h ,w \rangle_{\partial \refdom^m} &=(g\,f,w)_{\refdom} \,,\\
u_h &= \texttt{P}^{W_h}(\lambda_h),\qquad \textnormal{on }\Ff\, ,\\
\sum_{m=1}^M  \langle q_h  ,v \rangle_{\partial \refdom^m} &= \langle g_s\,h , v \rangle_{\partial \Omega_{\br,N}}  \, ,\end{align*}
holds for all $(w,v) \in W_h \times V_h(0)$, where the operator $\texttt{P}^{W_h}$ represents the projection onto the restriction of the space $W_h$ on the subdomain boundary. \RT{Note the last equation enforces continuity of the normal component of the flux across subdomain interfaces, as well as the Neumann condition on the Neumann boundaries}.


The next step is to eliminate the unknowns of $(u_h,q_h)$ to obtain a formulation only in terms of the variables  $\lambda_h$ defined on the subdomain interfaces. \RT{Invoking linearity and superposition, we consider} two subproblems for each subdomain $\refdom^m$: the first subproblem maps the function $f \in L^2(\phdom)$ to $u_h^f\vert_{\refdom^m} \in W_h^m(0)$, and the second subproblem maps $\eta\in V_h$ to $u_h^\eta \vert_{\refdom^m}\in W_h^m(\eta)$ as follows
\begin{alignat*}{2}
(\rho \bm G^m \nX u^f_h,\nX w)_{\refdom^m} - (g^m\,\wn u^f_h,w)_{\refdom^m}  &=(g^m\,f,w)_{\refdom^m} ,\quad &\forall w \in W_h^m(0),\\
(\rho \bm G^m \nX u^\eta_h,\nX w)_{\refdom^m} - (g^m\,\wn u^\eta_h,w)_{\refdom^m}  &= 0 ,\quad &\forall w \in W_h^m(0).
\end{alignat*}
Here $W^m_h = \{w \in \mcal{C}^0(\refdom^m),\,w|_T \in \mcal{P}^{p^m}(T), \forall T\in \mathcal{T}^m_h\}$ and $W^m_h(\eta) = \{w \in W^m_h : w = \texttt{P}^{W^m_h}(\eta),\; \textnormal{on } \partial \refdom^m\}$, where the operator $\texttt{P}^{W_h^m}$ represents the projection onto the restriction of $W_h^m$ on the subdomain boundary. Finally, we find the Lagrange multiplier $\lambda_h \in V_h(u_D)$ as a unique solution of the weak formulation 
\begin{equation*}
 (\rho \bm G \nX u^{\lambda_h}_h,\nX u_h^v)_{\refdom} -(g\,\wn u^{\lambda_h}_h,u^v_h)_{\refdom}  =(g\,f,u^v_h)_{\refdom} + \langle g_s\,h,v \rangle_{\partial\Omega_{\br,N}}
\end{equation*}
for $v\in V_h(0)$, and compute the numerical solution $u_h = u_h^f + u_h^{\lambda_h}$.

One of the main features of the MSCG method is the flexibility of mesh generation and \textit{h/p} adaptivity at the subdomain level, since the local subproblems are independent and are only coupled through the interface. Hence, the local subproblems may be statically condensed and expressed in terms of the variables defined on the subdomain faces. The idea is to relax the solution continuity at the subdomain interface and impose it back through a set of Lagrange multipliers $\lambda$, which correspond to the unique solution of the variational formulation on the mesh skeleton. Once the global solution is obtained, the solution for the local subproblems can be recovered independently for each subdomain. Another important advantage of the MSCG method is that the solution at the interfaces can be approximated with a different polynomial space than that of the subproblems, for instance using \RT{high-order Lagrange interpolation polynomials defined on the Chebyshev nodes}. This approach results in a drastic reduction of the size of the global system that would otherwise need to be solved had $(u_h,q_h)$ not been eliminated. Finally, continuity of the approximate solution across subdomains is enforced by imposing Dirichlet boundary conditions for each subdomain determined by the value of the global solution at the interfaces.

\subsubsection{Implementation}\label{sec:mscg_imp}

Let the total number of global (resp. local) degrees of freedom in the problem be given by $\mcal{N}_G$ (resp. $\mcal{N}_L$), and let $\mcal{N}^m_G$ (resp.  $\mcal{N}^m_L$) denote the number of global  (resp. local ) degrees of freedom associated to the $m$-th subdomain. Let the space $V_h$ be spanned by a set of global basis functions $\lbrace \varphi_i \rbrace_{i = 1}^{\mcal{N}_G}$. For each global face element $\mathsf{f}\in\Ff$, these global basis functions are chosen to be nodal Lagrange polynomials with the nodes placed at the Chebyshev points, as shown schematically in Fig. \ref{fig:globfun}. Then, we have $\lambda_h = \sum_{i=1}^{\mcal{N}_G} \Lambda_i \varphi_i(\refX)$, where ${\bf \Lambda} = (\Lambda_1,\ldots,\Lambda_{\mcal{N}_G})$ solves the global linear system
\begin{equation}\label{eq:globsystem}
 \mbb{K} {\bf \Lambda} = \mbb{F}\;,
\end{equation}
where $\lambda_h = \texttt{P}^{V_h}(u_D)$ is enforced on $\partial \Omega_{\br,D}$. The formation of the global stiffness matrix $\mbb{K}\in\mbb{C}^{\mcal{N}_G\times \mcal{N}_G}$ and load vector $\mbb{F}\in\mbb{C}^{ \mcal{N}_G}$ is detailed below.

Firstly, at the subdomain level we obtain $u_h^{f^m}$ and $u^{\varphi_i^m}_h \in W_h^m(\varphi_i^m),\; i = 1,\ldots, \mcal{N}_G^m$ by solving
\begin{subequations} \label{eq:subproblems_imp}
\begin{align}
(\rho \bm G^m \nX u^{f^m}_h,\nX w)_{\refdom^m} - (g^m\,\wn u^{f^m}_h,w)_{\refdom^m}  &=(g^m\,f,w)_{\refdom^m} \;,\label{eq:fproblem}\\
(\rho\bm G^m \nX u^{\varphi_i^m}_h,\nX w)_{\refdom^m} - (g^m\,\wn u^{\varphi_i^m}_h,w)_{\refdom^m} & = 0  \;,\label{eq:dirproblem}
\end{align}
\end{subequations}
for all $w \in W_h^m(0)$, where $\{\varphi^m_i\}_{i = 1}^{\mcal{N}_G^m}$ is the set of global basis functions that have non-zero support on $\partial\refdom^m$. A total of ${\mcal{N}_G^m}+1$ linear systems of dimension ${\mcal{N}_L^m}\times {\mcal{N}_L^m}$ are therefore solved for the $m$-th subdomain. The corresponding local stiffness matrix and load vector are computed as
\begin{subequations}\label{eq:locmatrix}
\begin{align}
 \mbb{K}_{ij}^m &= (\rho \bm G^m \nX u^{\varphi_j^m}_h,\nX u_h^{\varphi_i^m})_{\refdom^m} -(g^m\,\wn u^{\varphi_j^m}_h,u^{\varphi_i^m}_h)_{\refdom^m},\qquad 1\le i,j\le \mcal{N}_G^m, \label{eq:locmatrix1}\\
   \mbb{F}_{i}^m &= (g^m\,f,u_h^{\varphi_i^m})_{\refdom^m} + \langle g_s\,h,\varphi^m_i \rangle_{\partial \Omega_{\br,N} \cap \partial\refdom^m} ,\qquad 1\le i\le \mcal{N}_G^m. \label{eq:locmatrix2}
\end{align}
\end{subequations}
We then assemble $ \mbb{K}^m,\,\mbb{F}^m$ for $m = 1,\ldots,M$ on the global stiffness matrix $\mbb{K}$ and load vector $\mbb{F}$ through the standard finite element assembly procedure. After computing the global unknowns ${\bf \Lambda}$ from the linear system \eqref{eq:globsystem}, the solution on each subdomain ${\refdom^m}$ may be recovered as $u_h|_{\refdom^m} = u_h^{f^m} + \sum_{i=1}^{\mcal{N}_G^{m}} \Lambda_i^m u^{\varphi_i^m}_h$, where $\Lambda_i^m, 1 \le i \le \mcal{N}_G^{m}$ are the degrees of freedom of $\lambda_h$ that are nonzero on ${\refdom^m}$.

All in all, in the absence of geometry deformations (that is, $\mathfrak{G} = \mathbbm{1}$), only the local problem \eqref{eq:subproblems_imp}-\eqref{eq:locmatrix} for one representative of each subdomain class needs to be solved, since all the subdomains within the same class are identical. This feature is one of the key advantages of the MSCG method, which makes it particularly attractive for problems with periodic or repetitive structures.  Indeed, a judicious subdomain partition enables the simulation of complex photonic crystal structures just by solving a handful of small local subproblems and a single global linear system \eqref{eq:globsystem} defined only on the skeleton of the mesh.

However, for real applications the geometry of each subdomain within a class may not be exactly the same, see Fig. \ref{fig:motivation}, hence it would appear that this advantage is lost. In the next section, we introduce a reduced basis formulation that allows us to consider a more general parametrized problem that is, in fact, able to account for the small variations in the geometry that occur from subdomain to subdomain, and thus to partially retain the computational advantage of the MSCG method.

\subsection{3-D extension}
The formulation and implementation described above for the Helmholtz equation can be naturally extended to three-dimensional structures. However, there are some aspects of going from 2-D to 3-D that should be discussed.

The main obstacle is the increased computational cost that stems from considering linear systems that are larger and denser, that is going from 2-D to 3-D in solving the local subproblems \eqref{eq:subproblems_imp} and from 1-D to 2-D in solving \eqref{eq:globsystem} on the global face elements. In addition, since the Lagrange polynomials that approximate the global solution are now defined on surfaces, for each subdomain we have more global face elements and Dirichlet conditions per face. More precisely, using a tensorial grid on each global face squares the number of Dirichlet conditions. Even though sparse interpolation grids \cite{smolyak} can help mitigate this issue, the cost of \eqref{eq:dirproblem} augments significantly.

Finally, if the simulation of full three-dimensional photonic crystals is sought we can no longer reduce Maxwell's equations to Helmholtz, since the polarizations are coupled. Changing from Helmholtz to Maxwell is not straightforward from both the formulation and computation standpoints, since we now have to solve for six coupled three-dimensional fields instead of one, while ensuring the formulation is curl-conforming. To that end, a multiscale method for Maxwell's equations has been proposed and applied to 3-D optical fibers and 2.5-D photonic crystals in the dissertation \cite{saa2014simulation}, where the reader is referred for further details.

\section{Stochastic MSCG}\label{sec:sMSCG}

In this section, we introduce an approach for computing the statistics of quantities of interest, such as the transmission power of photonic devices. The main ingredients of our approach include: stochastic modeling to deal with geometry variability, a reduced order model \cite{barrault2004empirical,prud2002reliable,rozza2008reduced} to economize the simulation at the subdomain level, a variance reduction method \cite{giles2008multilevel,vidal2015model,vidal2016empirical} to accelerate the convergence of the statistical outputs, and an adjoint technique to compute derivatives.

In photonic crystals, the most widely used patterns are either dielectric circular rods in air or air holes in a dielectric slab \cite{joannopoulos2011photonic}. Ideally, all of these circular structures would be circles and have the same radius $R_0$. However, because of technological limitations at micro and nano scales, they are never perfectly round. For this work, we focus on photonic crystals with circular shapes, and model only geometry variations in the radial direction, albeit more geometries and complex deformations can be easily accommodated within our framework. The assumptions and mapping used for the examples in this article are thoroughly described in \ref{ap:map}.

\subsection{Reduced basis approximation}

A key component of our proposed approach is the use of a reduced basis (RB) approximation to efficiently compute the subdomain solutions for different rod geometries. For electromagnetism there has been considerable work to incorporate geometry as a parameter of the RB model \cite{chen2011seamless,hammerschmidt2016reconstruction,w2014reduced,hammerschmidt2016reduced}. \RT{Indeed, if the high-fidelity solutions lie in a low-dimensional manifold induced by the parametric dependence, then they can be accurately approximated using RB solutions at only a fraction of the cost, due to its rapid convergence. However, incorporating a reduced basis approximation in our multiscale framework is not straightforward, thus its construction is outlined below.}

The bilinear forms involved in the local problems \eqref{eq:subproblems_imp}-\eqref{eq:locmatrix} depend {\em non-affinely} on the parameter vector $\bm z$ since the deformation mapping $\mathfrak{G}$ is non-affine. To circumvent this issue, we use the empirical interpolation method (EIM) \cite{barrault2004empirical} at the subdomain level. For the $m$-th subdomain we approximate ${\bm G}^m$ and $g^m$ with the following affine expansions
\begin{equation}\label{eq:eim}
{\bm G}^m_Q(\refX,\bm z^m) = \sum_{q=1}^{Q_G} \sigma_q(\bz^m) \overline{\bm G}^m_q(\refX), \quad g_{Q}^m(\refX,\bz^m) = \sum_{q=1}^{Q_g} \varsigma_q(\bz^m) \overline{g}^m_q(\refX)\;.
\end{equation}
\RT{In our context, we compute for multiple values of $\bz^m$ the fields ${\bm G}^m$ and $g^m$ evaluated on the Gaussian quadrature points of the subdomain discretization, required to evaluate the volume integrals of the weak formulations \eqref{eq:subproblems_imp}-\eqref{eq:locmatrix}. Proper orthogonal decomposition (POD) \cite{berkooz1993proper} is then applied to attain the above expansions, keeping the modes necessary to retain at least $(1-10^{-8})\%$ of the system's energy, thus guaranteeing the approximation error is not dominated by the empirical interpolation. Naturally, the affine expansions may be reused for subdomains in the same class.} 

We then substitute \eqref{eq:eim} into \eqref{eq:subproblems_imp} and seek $u_h^{f^m} \in W_h^m(0)$ and $u_h^{\varphi_i^m} \in W_h^m(\varphi_i^m),\; i = 1,\ldots, \mcal{N}_G^m$ such that
\begin{subequations}
\begin{align} 
\sum_{q=1}^{Q_G}\sigma_q(\bz^m)(\rho \, \overline{\bm G}^m_q \nX u^{f^m}_h,\nX w)_{\refdom^m} - \sum_{q=1}^{Q_g}\varsigma_q(\bz^m)(\overline{g}^m_q  \,\wn u^{f^m}_h,w)_{\refdom^m}  &= \sum_{q=1}^{Q_g}\varsigma_q(\bz^m)(\overline{g}^m_q\,f,w)_{\refdom^m} ,\label{eq:eim1}\\
\sum_{q=1}^{Q_G}\sigma_q(\bz^m)(\rho \,\overline{\bm G}^m_q \nX u_h^{\varphi_i^m},\nX w)_{\refdom^m} - \sum_{q=1}^{Q_g}\varsigma_q(\bz^m)(\overline{g}^m_q \, \wn \,u_h^{\varphi_i^m},w)_{\refdom^m}  &= 0 .\label{eq:eim2}
\end{align}
\end{subequations}
for all $w \in W_h^m(0)$.

For subproblem \eqref{eq:eim2} assume that we are given orthonormalized basis functions $\zeta^{i}_n,\,1\le n\le N_{\max}$, such that $(\zeta^{i}_n,\zeta^{i}_{n'})_{W_h^m} = \RT{\int_{\refdom^m}\zeta^{i}_n\zeta^{i}_{n'} } =\delta_{n,n'}$. The associated hierarchical RB space is defined as $W^m_{i,N} = \mbox{span} \lbrace \zeta^{i}_n,\,1\le n \le N \rbrace$. In our case this spaces are constructed by computing a POD on a set of solutions. The main caveat of this formulation is that we typically have many Dirichlet conditions, which can impact the efficiency of the method since we must solve multiple RB systems. Instead, we propose to develop a single RB space $W^m_N$, where the basis functions $\zeta_n$ are computed with a weighted POD on a set of solutions for all possible Dirichlet boundary conditions. In other words, we treat the Dirichlet condition as an additional parameter, hence include the solutions for the various boundary conditions as independent snapshots. This strategy allows us to assemble and solve only one reduced basis system per subdomain group. For the first subproblem  \eqref{eq:eim1}, a similar process is followed if $f$ has spatial dependence. In such case, the RB needs to be constructed for all possible variations of the source on the subdomains. Nonetheless, if the source is constant, the RB for \eqref{eq:eim1} may be obtained using the standard RB procedure, and we denote the associated RB space as $\widetilde{W}^m_N(0)$.

We then apply a Galerkin projection to find $\mcal{N}_G^m+1$ RB solutions $u^{f^m}_N  \in  \widetilde{W}^m_N(0)$ and $u^{\varphi_i^m}_N \in W^m_N(\varphi_i) ,\; i = 1,\ldots, \mcal{N}_G^m$ satisfying, for all $(\widetilde{w},w) \in \widetilde{W}^m_N(0)\times W_N^m(0)$
\begin{subequations}\label{eq:rblocal}
\begin{align}
\sum_{q=1}^{Q_G}\sigma_q(\bz)(\rho \, \overline{\bm G}^m_q \nX u^{f^m}_N,\nX \widetilde{w})_{\refdom^m} - \sum_{q=1}^{Q_g}\varsigma_q(\bz)(\overline{g}^m_q\,\wn u^{f^m}_N,\widetilde{w} )_{\refdom^m}  &= \sum_{q=1}^{Q_g}\varsigma_q(\bz)(\overline{g}^m_q \,f,\widetilde{w} )_{\refdom^m},\\
 \sum_{q=1}^{Q_G}\sigma_q(\bz)(\rho \, \overline{\bm G}^m_q \nX u^{\varphi_i^m}_N,\nX w)_{\refdom^m} - \sum_{q=1}^{Q_g}\varsigma_q(\bz)(\overline{g}^m_q \,\wn u^{\varphi_i^m}_N,w)_{\refdom^m}  &= 0.
\end{align}
\end{subequations}
The local stiffness matrix $\mbb{K}^m_{ij}$ and load vector $\mbb{F}^m_{i}$ are then approximated as
\begin{subequations}\label{eq:rbelemglobal}
\begin{align}
 \mbb{K}^m_{ij} &\approx  \sum_{q=1}^{Q_G}\sigma_q(\bz^m)(\rho \, \overline{\bm G}^m_q \nX u^{\varphi^m_j}_N,\nX u_N^{\varphi^m_i})_{\refdom^m} - \sum_{q=1}^{Q_g}\varsigma_q(\bz^m)(\overline{g}^m_q\,\wn u^{\varphi^m_j}_N,u^{\varphi^m_i}_N\, )_{\refdom^m},\\
  \mbb{F}^m_{i} &\approx  \sum_{q=1}^{Q_g}\varsigma_q(\bz^m)(\overline{g}^m_q f,u^{\varphi^m_i}_N)_{\refdom^m} + \sum_{q=1}^{Q_g}\varsigma_q(\bz^m)\langle \overline{g}^m_q \,h,\varphi^m_i\rangle_{\refdomB{N}^m\cap \partial\refdom^m}.
\end{align}
\end{subequations}
for $1\le i,j\le \mcal{N}_G^m$. The affine parametric dependence enables an offline-online computational strategy where the parameter-independent instances are precomputed and stored beforehand, and the cost to obtain the approximate $u^{\varphi^m_i}_N$ for any parameter $\bm z^m$ and Dirichlet condition $\{\varphi^m_i\}_{i = 1}^{\mcal{N}_G^m}$ depends only on $Q_g,\, Q_G$ and $N$. 

\RT{It should be noted that, given the non-coercivity of the governing equation and the nature of the chosen projection, well-posedness of the reduced basis is not guaranteed. Projection techniques such as minimum-residuals \cite{maday2001blackbox} can be leveraged to regain well-posedness, at the expense of a higher computational cost. However, in the numerical experiments presented here the natural Galerkin projection described above the RB is well-posed, mainly due to the low frequencies of interest in photonic crystal design, thus the minimum-residuals formulation has not been used for this work.}

\subsection{Model and variance reduction}\label{ssec:mvr}
Our approach exploits the structure of the problem using the MSCG method and the RB for subdomains under a deformation mapping. We now study how these two ingredients can be combined to produce fast yet accurate estimates of the statistics of a quantity of interest. We propose a multilevel variance reduction method that exploits the statistical correlation between the different reduced basis approximations and the high-fidelity MSCG discretization to accelerate the convergence rate of Monte Carlo simulations. The multilevel variance reduction method provides efficient computation of the statistical outputs by shifting most of the computational burden from the high-fidelity MSCG approximation to the reduced basis approximations.  The methodology presented below is the direct extension of our previous work \cite{vidal2015model,vidal2016empirical} for the multiscale setting.

We first introduce a probability space $(\Upsilon,\Sigma,P)$, where $\Upsilon$ is the set of outcomes, $\Sigma$ is the $\sigma$-algebra of the subsets of $\Upsilon$, and $P$ is the probability measure. If $Z$ is a real random variable in $(\Upsilon,\Sigma,P)$ and $\upsilon$ a probability event, we denote its expectation by $E[Z] = \int_{\Upsilon} Z(\upsilon) d P (\upsilon)$. For an arbitrary subdomain $\refdom^m$, we shall consider functions $w\in L^2(\refdom^m \times \Upsilon)$ equipped with the following norm
\begin{equation*}
\norm{w}^2 = E\lb \int_{\mathcal{D}} |w( \refX,\cdot)|^2 d \refX\rb = \int_{\Upsilon} \int_{\mathcal{D}} |w( \refX,\upsilon)|^2 d \refX \,d P (\upsilon) .
\end{equation*}
We assume that, for a given subdomain $\refdom^m$, the parameters that define the geometric mapping $z^m_d(\upsilon)$ for $d= 1 \ldots, D$ are mutually independent random variables with zero mean. In addition, we assume that each of the $z^m_d(\upsilon)$ is bounded in the interval $\Gamma^m_d = [-\gamma^m_d, \gamma^m_d]$ with a uniformly bounded probability density function $\pi^m_d : \Gamma^m_d \to \mathbb{R}^+$. It thus follows that, with a slight overloading of notation, we can write $ \bm z^m = (z^m_1,\ldots,z^m_{D})$ and $\Gamma^m = \prod_{d = 1}^{D} \Gamma^m_d$. Hence, the entire stochastic space is given by $\bm\Gamma = \prod_{m = 1}^M \Gamma^{m}$ and the random variable as $\bm z = \prod_{m = 1}^M \bm z^{m}$. 

The solution $u$ of the original problem can be written as a function of $ \bm z \in \bm\Gamma$. Now let $s$ be a bounded functional. We introduce a random output $s$ defined as
\begin{equation*}
s( \bm z) = s(u(\cdot, \bm z)) .
\end{equation*}
We are interested in evaluating the expectation and variance of $s$ as
\begin{equation*}
 E[s] = \int_{\Gamma} s( \bm z) \pi( \bm z) d \bm z,\qquad \qquad 
V[s] = \int_{\Gamma} \lp E[s] - s( \bm z)\rp^2 \pi( \bm z) d \bm z,
\end{equation*}
where $\pi( \bm z) = \prod_{m=1}^M\prod_{d = 1}^{D} \pi_d^m(z^{m}_d)$. Since the exact output cannot be computed, we introduce the MSCG and RB outputs defined as $s_h(\bz),\,s_N(\bz)$ respectively. The high-fidelity output $s_h(\bz)$, computed by the MSCG method, is a very accurate approximation to $s(\bz)$ for any $\bz\in\Gamma$. Conversely, $s_N(\bz)$ denotes the output computed when we use an $N$-dimensional RB to evaluate \eqref{eq:rblocal}-\eqref{eq:rbelemglobal} on the subdomains with parametric dependence, as well as the high-fidelity CG formulation \eqref{eq:subproblems_imp}-\eqref{eq:locmatrix} for the remaining subdomains. \RT{Moreover, since the RB spaces provide a rapidly convergent approximation to the CG solution}, we expect a high statistical correlation between both outputs.

We now apply the above idea to compute an estimate of $E[s_h]$. To achieve this goal, we introduce 
\begin{equation*}
s_h^*(\bm z) = s_h(\bm z) + (E[s_{N}] - s_{N}(\bm z)),
\end{equation*}
where $s_{N}(\bm z)$ is the RB output for some $N \in [1, N_{\max}]$. Because $s_{N}(\bm z)$ generally approximates $s_h(\bm z)$ very well, the two outputs are highly correlated. The expectation may be recast as  
\begin{equation}\label{eq:2LE_exact}
E[s_h] = E[s_h^*] =  E[s_h - s_{N}] + E[s_{N}] .
\end{equation}
The underlying premise here is that the two expectation terms on the right hand side can be efficiently computed owing to variance reduction and model reduction: the first term requires a small number of samples because its variance is generally very small, while the second term is less expensive to evaluate because it only involves the RB output. The model and variance reduction (MVR) estimate of the expectation is computed through Monte Carlo simulations of the terms in \eqref{eq:2LE_exact} with $(I_0,\;I_1)$ i.i.d. samples, that is
\begin{equation}\label{eq:2LE_mvr}
E_{I_0,I_1}[s_h] = E_{I_0}[s_h - s_{N}] + E_{I_1}[s_{N}],
\end{equation}
which is an unbiased estimator of the expectation. 

The application of the CLT enables us to derive \apo estimate for the error in the expectation, that is
\begin{align}\label{eq:MVRerror}
\lim_{I_0 \to \infty}\lim_{I_1 \to \infty}& \mathrm{Pr} \left( \abs{E[s_h] - E_{I_0,I_1}[s_h]} \le \Delta_{I_0,I_1}^E \right) = \mathrm{erf}\lp\frac{\beta}{\sqrt{2}}\rp, \\
\Delta_{I_0,I_1}^E  &=  \beta \sqrt{\frac{V_{I_0}[s_h - s_{N}]}{I_0} + \frac{V_{I_1}[s_{N}]}{I_1}}\:,
\end{align}
where $\beta>0$ is the confidence level. For instance, in order to guarantee that $\abs{E[s_h] - E_{I_0,I_1}[s_h]}$ is bounded by a specified error tolerance $\epsilon_{\rm tol}$ with a high probability (say, greater than $0.95$), we need to take $\beta \ge 1.96$ according to the CLT. Similarly, the estimate of the variance is defined as
\begin{equation}\label{eq:2LV_mvr}
V_{I_0,I_1}[s_h] = E_{I_0}[\zeta_h - \zeta_{N}] + E_{I_1}[\zeta_{N}]\:,
\end{equation}
where the auxiliary variables are $\zeta_h  := \lp s_h - E_{I_0,I_1}[s_h] \rp^2$ and $\zeta_{N} := \lp s_{N} - E_{I_0,I_1}[s_h] \rp^2$. The error bound for the variance is defined analogously. The above derivations, the extension to the multilevel context and the optimal choice of $(I_0,\;I_1)$ and $N$ are thoroughly discussed and analyzed in \cite{vidal2015model}.

\subsection{Computation of gradients}
In the optimization context, the usage of first order optimization algorithms usually leads to accelerated convergence to the optimum value, as it guarantees a more efficient exploration of the design space. We now review how to obtain the gradients for the MSCG method. 

Following our assumptions, the long vector of parameters $\bz$ has at most dimension $MD$, since the subdomains not containing rods do not require any parameters. We now derive an adjoint approach to compute gradients in the multiscale context. The solution to the MSCG problem in matrix form is $\mathbf{u} = {\bf U}_{\varphi}\bm\Lambda + {\bf U}_f$, where ${\bf U}_{\varphi}\in\mbb{C}^{\mcal{N}_L \times \mcal{N}_G}$ is a matrix that contains the solutions to the Dirichlet subproblems and ${\bf U}_f$ is a vector of dimension $\mcal{N}_L$ that contains the solution to the source problem. Note that both ${\bf U}_{\varphi},\,{\bf U}_f$ are complex-valued, since the governing equations on PML subdomains are characterized by imaginary values of $\rho$ and $\kappa^2$ per \eqref{eq:pml}.

The derivatives of $s({\bf u},\bm z)$ can be recovered as
\begin{equation*}
 \frac{d s}{d \bm z} =  \frac{\partial s}{\partial \bm z} +  \frac{\partial s}{\partial \bf u}\frac{\partial \bf u}{\partial \bm z} = \frac{\partial s}{\partial \bm z} +  \frac{\partial s}{\partial \bf u} \lp \frac{\partial {\bf U}_{f}}{\partial \bz}+\frac{\partial \bf U_{\varphi}}{\partial \bm z} \bm \Lambda + {\bf U_{\varphi}}\frac{\partial \bm \Lambda}{\partial \bz}  \rp .
\end{equation*}
We invoke the adjoint technique to solve the last part, namely
\begin{equation*}
\frac{\partial s}{\partial \bf u} {\bf U_{\varphi}}\frac{\partial \bm \Lambda}{\partial \bz}  = \frac{\partial s}{\partial \bf u} {\bf U_{\varphi}} \mbb{K}^{-1} \lp \frac{\partial \mbb{F}}{\partial \bz}  - \frac{\partial \mbb{K}}{\partial \bz} \bm\Lambda \rp = \psi^\dagger\lp \frac{\partial \mbb{F}}{\partial \bz}  - \frac{\partial \mbb{K}}{\partial \bz} \bm\Lambda \rp,
\end{equation*}
where $^\dagger$ is the conjugate transpose operator, and the adjoint variable $\psi$ is the solution of
\begin{equation*}
 \mbb{K}^\dagger \psi = \lp \frac{\partial s}{\partial \bf u} {\bf U_{\varphi}} \rp^\dagger,
\end{equation*}
thus allowing us to compute the (at most) $MD$ derivatives with only an additional adjoint solution.

The derivatives of ${\bf U_{\varphi}},\, {\bf U_{f}},\,\mbb{K},\,\mbb{F}$ with respect to $\bz$ are sparse since the parameters only have local influence, and can be evaluated with \eqref{eq:subproblems_imp}-\eqref{eq:locmatrix} for the MSCG case and with \eqref{eq:rblocal}-\eqref{eq:rbelemglobal} for the RB case. An additional advantage of RB is the computation of these derivatives, since the affine dependence on the parameters that arises from the EIM greatly simplifies the computation.


\section{Numerical results}\label{sec:Res}
\RT{The MSCG algorithm introduced above is implemented in \MATLAB. To maximize computational efficiency, the routines to calculate the local stiffness matrix and load vector for CG \eqref{eq:locmatrix} and RB \eqref{eq:rbelemglobal} are fully vectorized. In order to present fair time comparisons, the same script is used to solve the entire multiscale problem, and the only difference is the routine used for local subproblems, either CG or RB. Thus, the global problem is assembled and solved identically for both CG and RB.  Time estimates are obtained averaging the wall time for 500 runs for each task using a single processor of a 512GB \Linux 18.04 machine with 32 AMD Opteron(tm) Processors 6320x15. For the results presented here, the multiple different local subproblems are solved in series, and parallelization is applied to the MVR computation of statistical moments given by \eqref{eq:2LE_mvr} and  \eqref{eq:2LV_mvr}.}

\subsection{Convergence test and cost analysis}
We start by studying the convergence of the MSCG method and compare it to that of the standard CG method for a simple problem. We consider problem \eqref{eq:helmholtz} in a unit square $\Omega = (0,1)^2$ for $\rho = \wn = 1$, where the Dirichlet condition and source term are chosen such that the exact solution is
\begin{equation*}
u(\bx) = x^2 + y^2 + \sin(k(x\cos\theta + y\sin\theta)) .
\end{equation*}
The exact solution is a plane wave of wavenumber $k$ propagating in the $\theta$ direction. The results reported here correspond to $\theta = \pi/4$ and $k = 6$.

For the CG discretization, we consider a homogeneous discretization of $2n^2$ triangular elements with uniform element size of $h = 1/n$. For the MSCG, we subdivide the domain into $q^2$ uniform squares, each of them with $2n^2/q^2$ triangles of uniform element size of $h = q/n$. To represent the solution, we consider both polynomial order $p=1,\,2$ for all subdomains; for the interfaces, besides considering a maximum polynomial order sufficient to capture the frequency of the problem, for instance 10, we also require for any boundary face that the local discretization is finer than the global discretization in order to capture the boundary conditions prescribed by the Lagrange polynomial. Hence we take $p^{\mathsf{f}} = \min \lbrace 10,np/q \rbrace$, and discretize each global face with a single high-order element.

We compute the $L^2(\Omega)$ errors for both methods for $n = 8,16,32,64,128$. Since the global polynomial order suffices for the frequency of interest, the error is dominated by the subproblems, where both methods have the same discretization, hence both the \RT{CG and MSCG errors for different $q$ values coincide}, and yield the expected convergence rate of $\mcal{O}(h^{p+1})$ for smooth solutions, see Table \ref{tab:p1}. 

\begin{table}[h!]
\footnotesize
\begin{center}
  \renewcommand{\arraystretch}{1.2}
\begin{tabular}{cccccc}
   & \multicolumn{2}{c}{$p=1$} &  &\multicolumn{2}{c}{$p=2$} \\
$n$  & $\norm{u-u_h}_{L^2}$ & Order & &$\norm{u-u_h}_{L^2}$ & Order \\
8  & 2.80e-2 & -- & &6.23e-4 &  -- \\
16 & 7.63e-3 & 1.86 & &7.55e-5 & 3.04 \\
32 & 1.96e-3 & 1.97 & &9.38e-6 & 3.01 \\
64 & 4.92e-4 & 1.99 & &1.17e-6 & 3.00 \\
128& 1.23e-4 & 2.00 & &1.46e-7& 3.00 \\
\end{tabular}
\end{center}
\caption{$L^2(\Omega)$ error and convergence rate for $p = 1,\,2$ as a function of $n$.}\label{tab:p1}
\end{table}

Another relevant study is the comparison of degrees of freedom for both methods. A uniform triangular mesh of $2n^2$ elements of order $p$ renders $\mcal{N}_{\rm{CG}} = (np+1)^2$ high-order nodes. Conversely, for each of the $q^2$ subdomains of the MSCG we have $2(n/q)^2$ elements and $\mcal{N}_L = (np/q+1)^2$ degrees of freedom. Furthermore, if we assume a homogeneous polynomial order $p^\mathsf{f}$ for the Lagrange polynomial at the interfaces, it can be shown that $\mcal{N}_G = (q+1)(2qp^\mathsf{f}-q+1)$. If we adopt $p^{\mathsf{f}} = \min \lbrace 10,np/q \rbrace$, which suffices to represent the solution accurately for this case, we obtain the degrees of freedom presented in Table \ref{tab:p2dof} for $p=2$.

\begin{table}[h!]
\footnotesize
\begin{center}
  \renewcommand{\arraystretch}{1.2}
\begin{tabular}{cccccccc}
   & \multicolumn{6}{c}{ MSCG} &  CG  \\
  & \multicolumn{2}{c}{$q=2$}& \multicolumn{2}{c}{$q=4$} & \multicolumn{2}{c}{$q=8$} & \\ 
$n$& $\mcal{N}_L$ & $\mcal{N}_G$ & $\mcal{N}_L$ & $\mcal{N}_G$ & $\mcal{N}_L$ & $\mcal{N}_G$ &$\mcal{N}_{\rm{CG}}$\\
8  & 81 & 93 & 25 & 145 & 9 &225 & 289\\
16 & 289 & 117 & 81 & 305 & 25 & 513 & 1089\\
32 & 1089 & 117 & 289 & 385 & 81 & 1089 & 4225\\
64 & 4225 & 117 & 1089 & 385 & 289 & 1377 & 16641\\
128& 16641 & 117 & 4225 & 385 & 1089 & 1377 & 66049\\
\end{tabular}
\end{center}
\caption{Degrees of freedom for MSCG and CG for $p = 2$ as a function of $n$.}\label{tab:p2dof}
\end{table}
This analysis shows that, for this problem where there is only one subdomain type, the best strategy is to use more subdivisions per direction as we refine the mesh. Indeed, a multiscale configuration is deemed optimal whenever the number degrees of freedom for both the local and global problems are similar, thus the computational burden is divided evenly. Furthermore, notice that the benefit of the MSCG as opposed to regular CG becomes apparent as more elements are used in the discretization. The results presented here correspond to the very simple case where only one local computation needs to be performed, therefore if there is more than one subdomain type one needs to account for the multiple unique local subproblems. 

\subsection{Waveguide T-splitter simulation}

We now focus on the simulation of waveguides in photonic crystals. Photonic crystals are assembled by combining different materials, giving rise to periodic nanostructures that exhibit the bandgap phenomenon. That is, there exists broad bands of frequencies for which wave propagation through the crystal is disallowed. Common examples of such structures are obtained by uniformly placing dielectric rods in air or drilling regular patterns of holes in a dielectric slab. The bandgap response of these structures has been studied extensively \cite{joannopoulos2011photonic}. Waveguiding arises as an application of photonic crystals whenever symmetry of the lattice is broken, for instance when a line of rods is removed from the crystal. In this scenario, if the crystal is illuminated with a wave whose frequency is in the gap, the wave will only be allowed to travel along the defect and will decay exponentially away from the defect. The MSCG presented in this paper is therefore an attractive candidate to simulate photonic crystal applications, since by solving a small number of subproblems it enables the simulation of large lattice structures by exploiting the repeated patterns. 

\begin{figure}[h!]
 \centering
 \includegraphics[scale = .45]{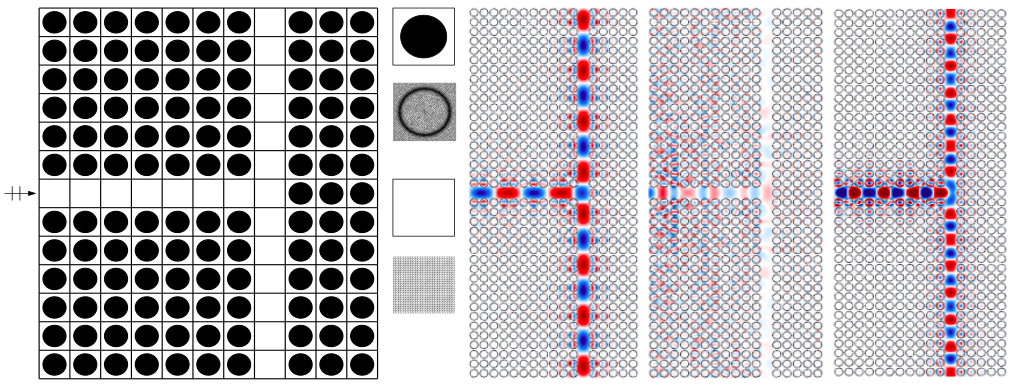}
 \caption{MSCG simulation of a TM waveguide splitter of GaAs ($\varepsilon = 11.4$) rods of radii $R_0 = 0.4a$ in air. Computational domain with subdomain decomposition (excluding PMLs) and meshes (left). Numerical simulation at frequencies $\overline{\omega}= 0.39$, $\overline{\omega}= 0.46$ and $\overline{\omega}= 0.53$ (left to right).}\label{fig:mscg_splitter}
\end{figure}

The first example is a waveguide T-splitter consisting of Gallenium Arsenide rods ($\varepsilon = 11.4$) in air of radii $R_0 = 0.4a$, where $a$ is the periodicity of the crystal. For the TM polarization, this structure presents a first bandgap for $\overline{\omega}\in (0.36,0.40)$ and a second bandgap for $\overline{\omega}\in (0.52,0.55)$. In order to numerically simulate the splitter with the MSCG method, we first have to identify the subdomains in which to split the computational domain. The subdomains should be invariant to rotation and translation to ensure that only a small number of subproblems are solved. The size of the subdomain along with the frequencies of interest determine the polynomial order chosen to approximate the solution at the global interfaces. Conversely, the discretization needed at the subdomain level is governed by the details of the geometry. Furthermore, the MSCG approach offers the possibility of using different polynomial approximation order for different subdomains, which can result in greater efficiency.

For the waveguide splitter shown in Fig. \ref{fig:mscg_splitter} (left), we choose three different classes of subdomains:  (1) the rod subdomain, where we use $p=2$, curved elements and adaptive mesh size to approximate the fine detail of the rod curvature; (2) the air subdomain, where we use an homogeneous mesh of straight-sided elements of order $p=2$; and (3) the PML subdomain, with the same discretization as type 2. Although not represented in Fig. \ref{fig:mscg_splitter}, three PML subdomains surround the computational domain in each direction to ensure outgoing waves are properly attenuated. \RT{The subdomain discretization is chosen through a grid convergence study on progressively refined meshes, until the relative error of the optical intensity at both output ports is below 0.1\%, compared to an extremely fine subdomain discretization.} For each global interface, we choose a single element of polynomial order $p^\mathsf{f} = 10$, which gives a sufficient resolution of more than 20 points per wavelength. The degree of freedom count is reported in Table \ref{tab:dofsplitter}. Note that solving the exact same problem with regular CG would suppose solving a linear system of more than 6M dof as opposed to a system of 26K dof, thus for structured geometries the MSCG is a very competitive approach.

\begin{table}[h!]
\footnotesize
\begin{center}
  \renewcommand{\arraystretch}{1.2}
\begin{tabular}{cccccc}
subdomain type  & $\#$  elements/subdomain &  order & $\#$  dof per subdomain & {$\#$ subdomains} &  $\mcal{N}_L$   \\
1 & 3312 & 2 & 7K & 613 & 4M\\
2 & 1352 & 2 & 3K & 50 & 150K\\
3 & 1352 & 2 & 3K & 663 & 1.8M\\
\hline
global faces & $\#$  elements/face &  order &$\#$ dof per face & {$\#$ faces} &  $\mcal{N}_G$   \\
 & 1 & 10 & 11 & 2729 & 26K\\
\end{tabular}
\end{center}
\caption{Degrees of freedom for MSCG for the waveguide T-splitter separated by subdomain types.}\label{tab:dofsplitter}
\end{table}

In Fig. \ref{fig:mscg_splitter} (right) we show the amplitude field $E_z$ for a frequency in the first bandgap, a frequency between the two bandgaps and a frequency in the second bandgap. An attractive feature of waveguiding with photonic crystals are the low losses that occur even for sharp bends, thus enabling the efficient manipulation of electromagnetic waves.

\subsection{Waveguide Z-bend optimization}

\begin{figure}[h!]
\centering
\subcaptionbox{Left: Silicon slab with Z-bend defect on a triangular lattice. Right: Computational domain with different subdomains classes highlighted. The red subdomains correspond to holes in the optimization area.\label{fig:slabDomain}}
[15cm]{\includegraphics[scale = 0.45]{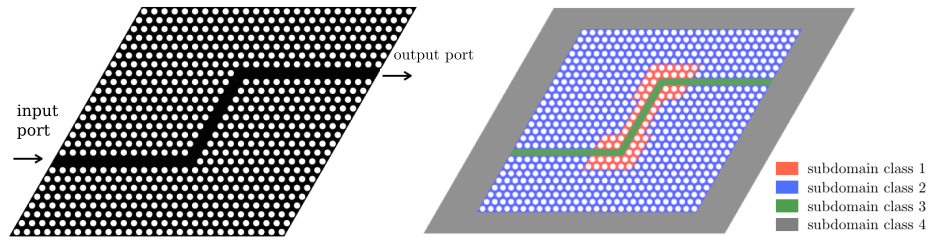}}
\hfill\subcaptionbox{Quantity of interest vs frequency, with amplitude field shown for $\overline{\omega}\in\lbrace 0.2655,0.285,0.305\rbrace$. \label{fig:slabSolutions}}
[15cm]{\includegraphics[scale = 0.35]{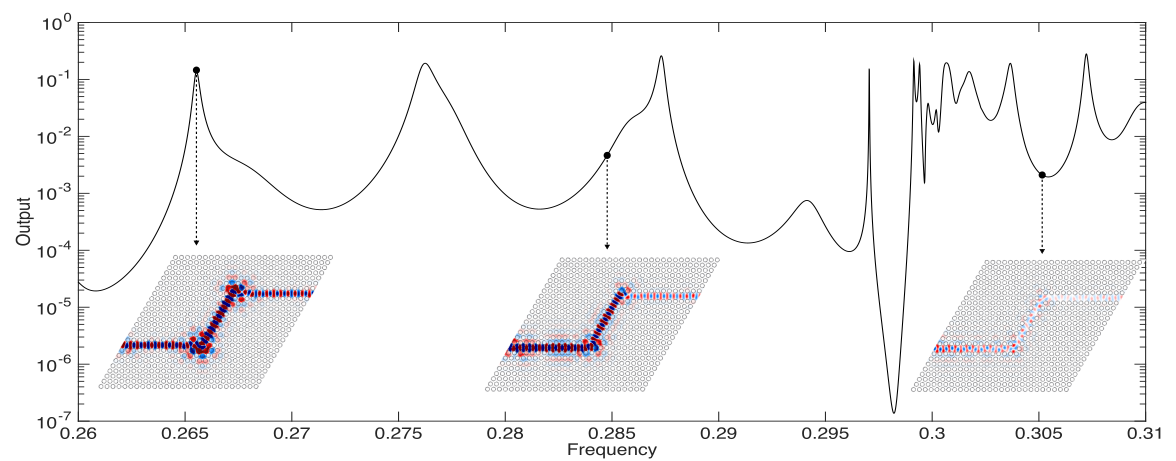}}
\caption{Schematics, computational domain and solutions for the silicon slab wave propagation problem.}
\end{figure}

Here we consider the simulation and design of a photonic crystal consisting of a silicon ($\varepsilon = 12.1$) slab with drilled air holes of radius $R_0 = 275a/800$ forming a lattice with triangular symmetry. This structure presents a broad bandgap for $\overline{\omega}\in (0.26,0.34)$ for TE waves, see \cite{frandsen2004broadband}. A waveguide is generated by opening a Z-shaped defect  and illuminating the crystal with a wave impinging at the left input port. The quantity of interest is the intensity of the optical power at the right output port in the $x$-direction, namely
\begin{equation*}
s_h = \frac{1}{2\omega}\sum_{j=1}^J \Bigl\lvert \int_{\Omega^j_{out}} \varepsilon^{-1} \,{\bf e}_x\cdot\Re \lb \mathrm{i} u \nabla u^{\dagger}\rb \Bigr\rvert,
\end{equation*}
where the subdomains at the output port $\Omega_{out}$ comprise the line defect and one hole subdomain above and below. The schematics of the slab are shown in Fig. \ref{fig:slabDomain} (left), where we have ensured that the bends are sufficiently separated such that they do not interact. The objective here is to apply the methodology introduced above to design the bending region such that transmission is enhanced for certain frequencies of interest within the bandgap of the regular lattice \cite{frandsen2004broadband}, see Fig. \ref{fig:slabSolutions}.


The different subdomain classes identified are shown in Fig. \ref{fig:slabDomain} (right). The subdomain classes 1 and 2 correspond both to the hole subdomain, but are treated independently since class 1 will be subject to shape optimization, whereas the subdomains from class 2 remain unaltered. The other two subdomain types correspond to the line defect and the PML respectively. The dof count for the different types and the global problem is collected in Table \ref{tab:dofslab}.

\begin{table}[h!]
\footnotesize
\begin{center}
  \renewcommand{\arraystretch}{1.2}
\begin{tabular}{cccccc}
subdomain class  & $\#$  elements/subdomain &  order & $\#$  dof per subdomain & {$\#$ subdomains} &  $\mcal{N}_L$   \\
1 & 9196 & 2 & 19K & 58 & 1M\\
2 & 9196 & 2 & 19K & 607 & 11M\\
3 & 4050 & 2 & 8K & 37 & 0.3M\\
4 & 4050 & 2 & 8K & 482 & 4M\\
\hline
global faces & $\#$  elements/face &  order &$\#$ dof per face & {$\#$ faces} &  $\mcal{N}_G$   \\
 & 1 & 10 & 11 & 2437 & 23K\\
\end{tabular}
\end{center}
\caption{Degrees of freedom for MSCG for the waveguide Z-bend separated by subdomain types.}\label{tab:dofslab}
\end{table}

%

In the next sections, we present optimization results for three different scenarios: single frequency, range of frequencies and range of frequencies combined with the presence of geometric imperfections in the holes. In all cases, the optimization variables are the radii of the holes in subdomain class 1. Hence, we have $M_1 = 58$ design parameters $\theta^{m}\in \Theta^m$ defined as $r_0^m = R_0(1+\theta^m),\, 1\le m \le M_1$. Box constraints are prescribed for all design parameters, and are set as $\Theta^m = [-0.127,0.047]a$ such that $r^m \in [0.3,0.36]a$, hence the design region is $\bm\Theta:=\Theta^1\times\ldots\times\Theta^{M_1}$. 

The optimizations are performed using the \texttt{nlopt} \cite{nlopt}  optimization package. \RT{For each objective function, we first apply the multi-level single-linkage global optimization algorithm, and then choose the computed optimum as the starting point for a local gradient-based optimization using the preconditioned truncated Newton, in order to perform a better local exploration, see references in \cite{nlopt}. Relative tolerances for both the objective function value and the design parameters are employed as termination criteria, and are set to $10^{-2}$ for the global optimization and $10^{-4}$ for the local optimization. Finally, for each objective function we compute multiple global-local optimization cycles with different initial guesses, in order to avoid getting stuck in local minima. The results shown here correspond to the optimization cycle that achieved the lowest objective function value.}

\subsubsection{Single frequency radii optimization}
The first step towards efficient optimization is to develop a RB for subdomain class 1, which reduces to finding a RB for the Dirichlet subproblem \eqref{eq:dirproblem} as we have no source. The mapping for subdomain class 1 is described in \ref{ap:map}, whereas the mapping for the other subdomain classes is the identity. The application of the EIM on the three components of $\bm G$ (represents a second-order symmetric tensor) and the Jacobian $g$ enables the affine parametrization of the mapping, for a total of $Q_G=6$ and $Q_g=4$ interpolation elements. We then compute snapshots to \eqref{eq:dirproblem} for the different Dirichlet conditions, where the parameters are $\bm\theta$ and the frequency $\overline{\omega}\in [0.3045,\,0.3055]$ --- thus the RB can be reused for the next example. Finally, the RB is constructed by compressing the snapshots for all the boundary conditions to form a single POD basis.

We now proceed to optimizing the radii of the holes for a single frequency $\overline{\omega}= 0.305$, for which intensity at the outlet is poor. For subdomains 2, 3 and 4 we use CG to compute \eqref{eq:dirproblem}-\eqref{eq:locmatrix}  in 0.42, 0.15 and 0.30 seconds respectively.  Thus for each new $\bm \theta$ we only need to solve for the 58 subdomains from class 1.  The deterministic optimization problem reads 
\begin{equation}\label{eq:DO}
s_N(\bm\theta^*) = \max_{\substack{\bm\theta \in \bm\Theta}} \quad s_N\;,
\end{equation}
where the expensive CG local solution \eqref{eq:dirproblem}-\eqref{eq:locmatrix} is substituted for its inexpensive RB counterpart \eqref{eq:rblocal}-\eqref{eq:rbelemglobal} to accelerate the computations. \RT{We use a reduced basis of size $N=100$, rendering a wall time for all 58 suproblems of 0.12 s, as opposed to 25.60 s required by the high-fidelity CG method. The global system is finally solved in 0.47 s. Overall, the reduced basis approach delivers the value of $s_N$ in just 1.45 s, that is 18 times faster than computing $s_h$ with CG on all subdomains, and without compromising accuracy as shown below.}

Convergence is achieved after 1200 iterations. A frequency sweep for $s_N(\bm\theta^*)$ is shown in Fig. \ref{fig:expectation} (dashed red), showing a much higher transmission than the un-optimized version, see Fig. \ref{fig:slabSolutions}. \RT{Lastly, accuracy of RB is verified by re-optimizing \eqref{eq:DO} with $s_h$ instead of $s_N$, setting the maximizer $\bm\theta^*$ as initial guess.  The new optimum $\bm\theta_h^*$ is attained only after 3 iterations, since $s_N$ is a faithful representation of $s_h$. The relative error of the output for the same maximizer, $\abs{s_h(\bm\theta^*)-s_N(\bm\theta^*)}/s_h(\bm\theta^*) = 5.8e-4$, as well as the relative error of the maximizer, $\norm{\bm\theta_h^* - \bm\theta^*}_2/\norm{\bm\theta_h^*}_2 = 8.9e-4$, are below the acceptable engineering standards, therefore concluding the RB effectively approximates the high-fidelity CG solutions of the local subproblems.}

\subsubsection{Frequency-range radii optimization}

The next step is to seek a design that is robust with respect to a range of frequencies, which is desirable when the optimum presents a sharp peak with a rapid decay away from the single frequency optimum. This optimization problem can be recast in a stochastic optimization framework as
\begin{equation}\label{eq:SOfreq}
\widehat{s}_N(\widehat{\bm\theta}) = \max_{\substack{\bm\theta \in \bm\Theta}} \quad E_\omega [s_N] - \gamma \sqrt{V_\omega[s_N]}
\end{equation}
where $\gamma$ controls the weight assigned to the variance minimization. Since the stochastic space is unidimensional, the output statistics can be computed with simple Gauss-Legendre quadrature, reducing \eqref{eq:SOfreq} to the weighted evaluation of the output $s_N$ at selected frequencies given by Legendre points. \RT{We choose a 20-point quadrature that integrates exactly polynomials up to degree 39, selected as the coarsest quadrature that renders a mean relative error on \eqref{eq:SOfreq} of $<0.1\%$ for multiple $\bm\theta$ values, compared to the results of a 100-point quadrature.} We reuse the RB constructed in the previous section for subdomain class 1, again with 100 modes, and precompute \eqref{eq:dirproblem}-\eqref{eq:locmatrix} for subdomain classes 2, 3 and 4 at the selected Legendre frequencies in the range $\overline{\omega}= [0.3045,0.3055]$ to expedite the optimization process. A frequency sweep for $s_h(\widehat{\bm\theta})$ (dashed green) with $\gamma = 1$ is shown in Fig. \ref{fig:expectation}. The radii configuration $\widehat{\bm\theta}$ is indeed less sensitive to frequency variations in the prescribed range, albeit attaining a lower output intensity at  $\overline{\omega}= 0.305$ compared to the single-frequency optimum ${\bm\theta}^*$ (dashed red). 

\RT{Similarly as before, we can assert the RB is accurate by inspecting both the relative error on the optimum if the objective function is evaluated with CG ($\widehat{s}_h$) instead of RB ($\widehat{s}_N$) for the subproblems belonging to class 1, namely $\abs{\widehat{s}_h(\widehat{\bm\theta})-\widehat{s}_N(\widehat{\bm\theta})}/\widehat{s}_h(\widehat{\bm\theta}) = 1.8e-4$, and the relative error of the optimum $\widehat{\bm \theta}_h$ when maximizing $\widehat{s}_h$ initialized with $\widehat{\bm\theta}$, that is $\norm{\widehat{\bm\theta}_h - \widehat{\bm\theta}}_2/\norm{\widehat{\bm\theta}_h}_2 = 7e-4$.}

\subsubsection{Frequency-range robust radii optimization}

Finally, we analyze the robustness of the solutions with respect to geometry errors on a range of frequencies. The geometry errors, which will be considered for subdomains in classes 1 and 2, are given by \eqref{eq:KL}, where we append the constant deformation $\theta^m$ for the optimization, namely  $r_0^m \mapsto R_0^m\lb 1 + \theta^m\rb + \delta R_0^m$. For the results below, we select $D = 11$, $\sigma = 0.02$, $L_c = 1/16$  and $\bm z^m \in [-\sqrt{3},\sqrt{3}]^{D}$ in \eqref{eq:KL}, for which the 95\% confidence interval gives $\abs{\delta R_0^m}/R_0^m < 3\%$ if $\theta = 0$. Note we encapsulate both the optimization and the stochastic parameters in the same expression, and thus the same RB for both classes of subdomains can be used by simply setting $\theta^m=0$ for the subdomains in class 2.

In order to achieve robust designs accounting for variation in both the frequency and the geometry parameters, we formulate the following stochastic optimization problem 
\begin{equation}\label{eq:SOgeom}
\widetilde{s}_h(\widetilde{\bm\theta}) = \max_{\substack{\bm\theta \in \bm\Theta}} \quad E_{\omega,\mathcal{G}} [s_h] - \gamma \sqrt{V_{\omega,\mathcal{G}}[s_h]}.
\end{equation}
The stochastic dimension of the problem under consideration is large, since we have 11 geometry parameters per subdomain, and there are 665 subdomains of classes 1 and 2 combined. Hence, we shall resort to the MVR method described in Section \ref{ssec:mvr} to evaluate the objective function \eqref{eq:SOgeom}. We first develop a RB for the Dirichlet subproblem combining the EIM with POD for the non-homogeneous radius variation. We set $\theta^m \in [-0.127,0.047]a$ as before for the optimization, which in this case leads to $Q_G = 61$ and $Q_g=23$ elements in the interpolation basis. We again use $N=100$ for the RB model, but since have considered 11 additional parameters to represent the geometry we can expect this RB to be less accurate than the one developed for the homogeneous radii variation. The main difference here is that the RB is seen as a surrogate that correlates with the high-fidelity model, not as a substitute, thus a coarser basis still casts excellent results. 

\RT{Subdomains 3 and 4 are solved with CG in 0.15 and 0.30 seconds respectively, and the global system in 0.47 s. Unfortunately, the CG computation for subdomain class 2 can no longer be precomputed and reused due to the geometry errors, thus requiring 665 high-fidelity CG linear system solutions for subdomains 1 and 2 combined. This calculation takes 301.3 s in our implementation, for a total MSCG wall time of 302.2 s for a single tuple $(\bm\theta,\bm z)$. Additionally, computing the output statistics at each optimization step typically demands multiple output evaluations for different $\bm z$, thereby greatly hindering the computational feasibility of the optimization. To circumvent this issue we employ the RB on subdomains 1 and 2, reducing the wall time to 3.22 s for a total MSCG time of 4.16 s. The MVR method combines the fact that $s_N$ is more efficiently evaluated than $s_h$ (a factor of 70), with high variance reduction between outputs $V(s_h-s_N)/V(s_h) > 50$, therefore conveniently allocating many more samples to the RB output than to the variance reduction term ($I_1 = 15K$ and $I_0 = 355$ for $\Delta^E_{I_0,I_1} = 10^{-3}$ at the optimum) in order balance the computational cost of \eqref{eq:2LE_mvr}.}

\begin{figure}[h!]
\centering
\subcaptionbox{Output $s_h(\bm\theta)$ (dashed) and expected value of output with respect to geometry with 95\% CI $E_{\mcal{G}}[s_h(\bm\theta)] \pm \Delta^E$ (solid) for the three optimized designs in the frequency range of study.\label{fig:expectation}}
[7.5cm]{\includegraphics[scale = 0.35]{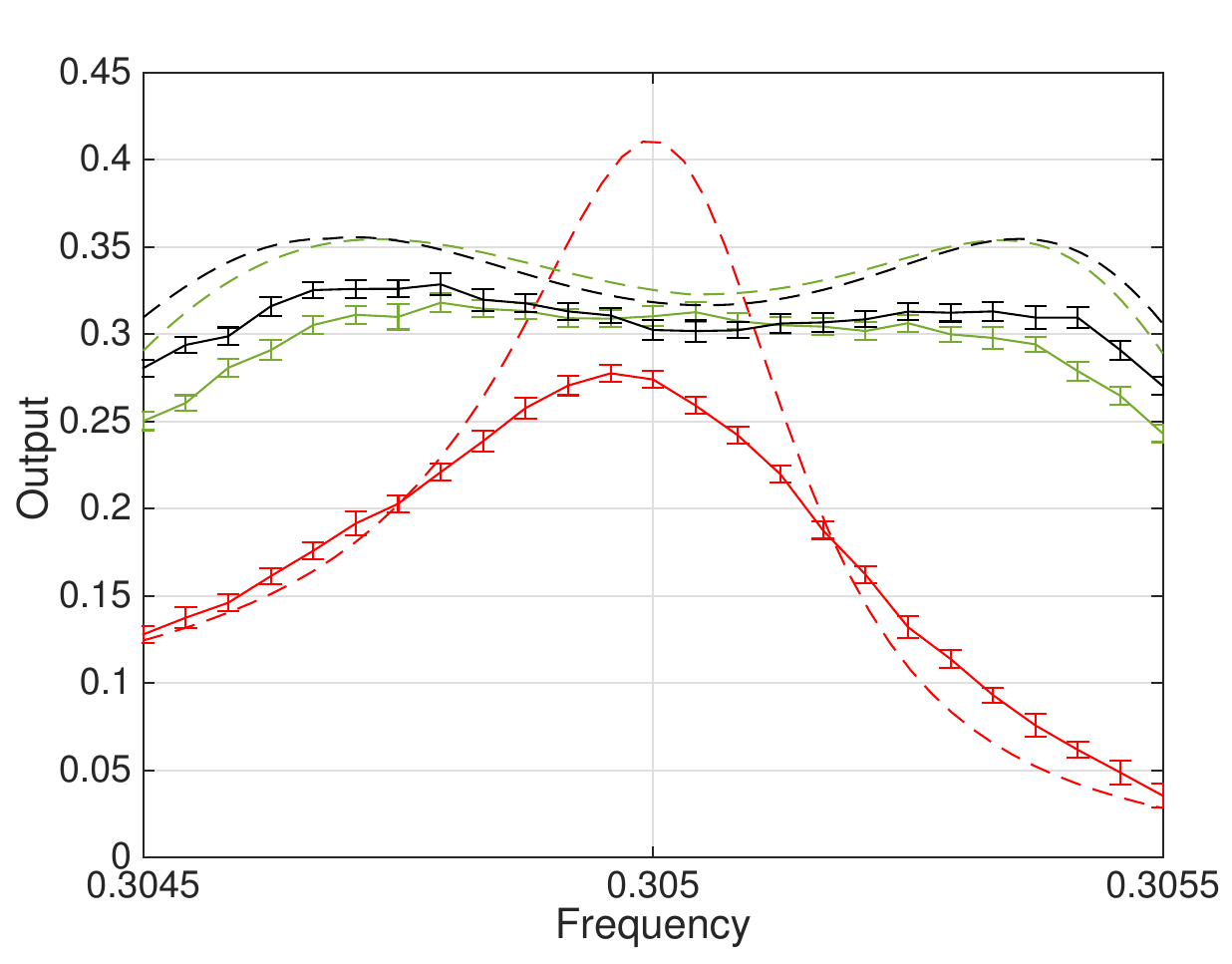}}
\hfill\subcaptionbox{Variance of output with respect to geometry with 95\% CI $V_{\mcal{G}}[s_h(\bm\theta)] \pm \Delta^V$ for the three optimized designs in the frequency range of study. \label{fig:variance}}
[7.5cm]{\includegraphics[scale = 0.35]{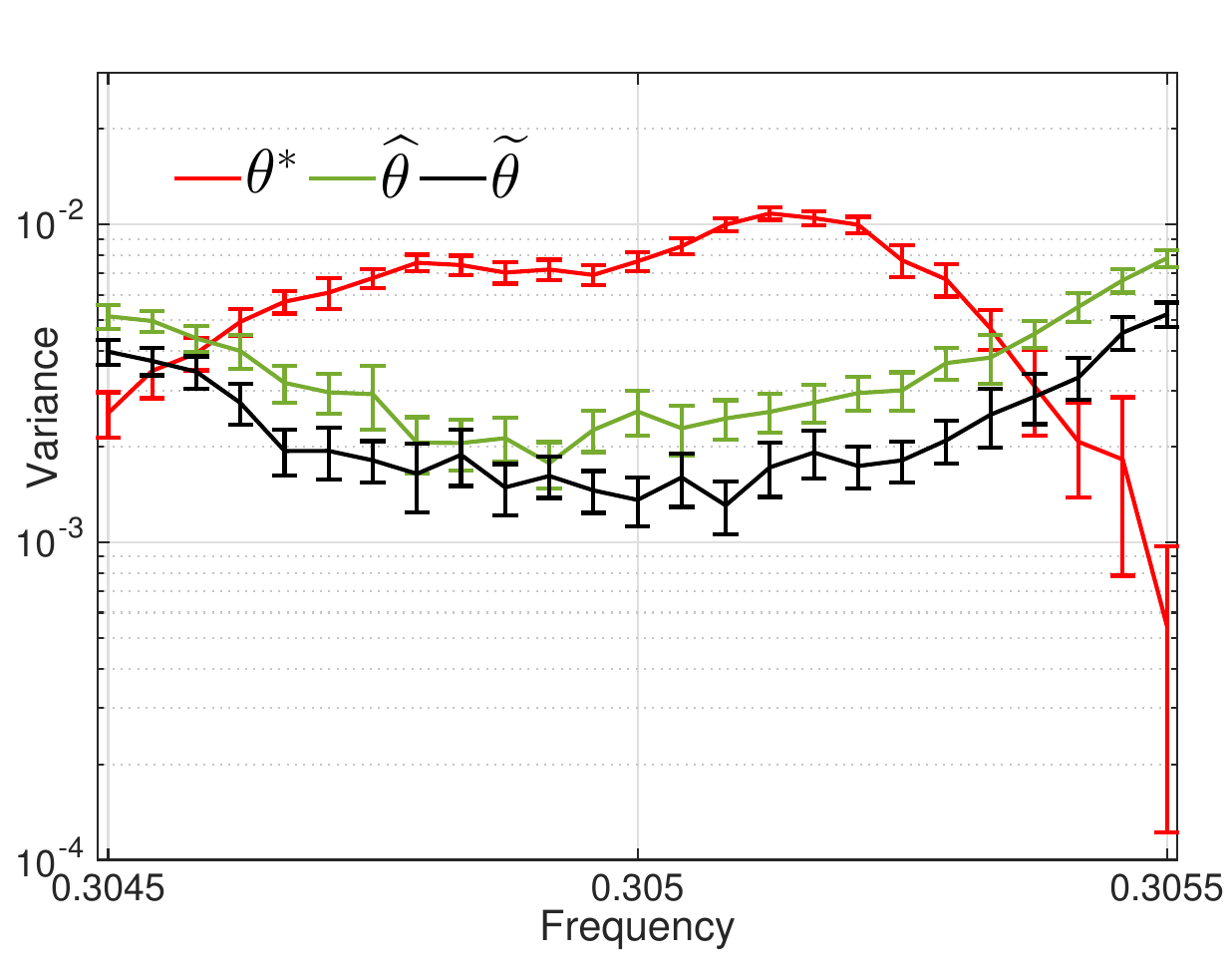}}
\hfill\subcaptionbox{Single frequency $\bm\theta^*$ (left), range of frequencies $\widehat{\bm\theta}$ (middle) and robust design $\widetilde{\bm\theta}$ (right) optima. \label{fig:designs}}
[15cm]{\includegraphics[scale = 0.35]{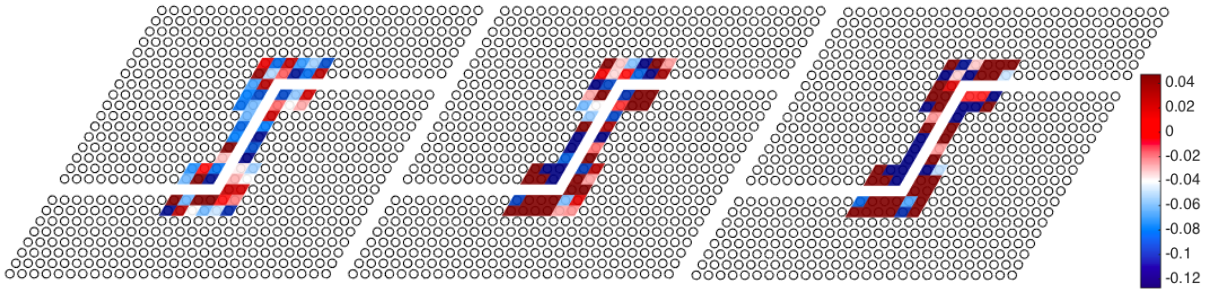}}
\caption{Results for the waveguide Z-bend optimization problem.}
\end{figure}

In Fig. \ref{fig:expectation} we show in dashed line, for a set of discrete frequencies in the interval $\overline{\omega}= [0.3045,0.3055]$, the transmission $s_h$ for the several configurations considered: single frequency optimization $\bm\theta^*$, frequency-range $\widehat{\bm\theta}$ and frequency-range robust design $\widetilde{\bm\theta}$. To assess the robustness with respect to geometry, for each design we introduce geometry variations on the subdomain classes 1 and 2, and evaluate the expected value of the output and its 95\% confidence interval as $E_{\mcal{G}}[s_h(\bm\theta)] \pm \Delta^E$ in Fig. \ref{fig:expectation}, as well as the variance $V_{\mcal{G}}[s_h(\bm\theta)] \pm \Delta^V$ in Fig. \ref{fig:variance}, both represented in solid line. The single frequency optimum  $\bm\theta^*$ produces the highest transmission for  $\overline{\omega}= 0.305$, but it degrades significantly in the presence of geometric errors. Conversely, the other optima, despite achieving lower peak outlet intensities, maintain a consistent performance for all the frequencies in the range. Moreover, the robust optimum  $\widetilde{\bm\theta}$ outperforms the other two optima both in expected value and in variance, see Fig. \ref{fig:variance}. These results show the importance of accounting for geometry variations in the objective function if robust designs with respect to geometry are sought. Finally, the optima configurations for the three scenarios considered are depicted in Fig. \ref{fig:designs}, expressed as a variation on the nominal radius.

\section{Conclusions}\label{sec:conclusions}

We have presented a multiscale continuous Galerkin method for computing wave propagation phenomena in photonic crystals. The method relies on partitioning the computational domain into subdomains, and computing the local solution in the subdomains as a function of the Lagrange multipliers at the interfaces. Thus, only a linear system for the global Lagrange multipliers is solved. The solution in the interior of the subdomains is recovered after the global problem has been solved. Furthermore, the MSCG method is especially advantageous for problems that exhibit repeated patterns, since the local subproblems only need to be solved once for each type of subdomain. In addition, we propose the construction of a reduced order model for the local subproblems that allows us to consider geometric variation within the subdomains, leveraging a reference domain formulation of the governing equation. We then use the model and variance reduction method to compute statistical outputs of stochastic wave propagation problems. Finally, we verify the convergence rate of MSCG, carry out a deterministic simulation of a waveguide splitter to illustrate the advantages of MSCG in photonic crystal, and a robust optimization of a photonic slab to demonstrate the performance of the MSCG method to attain robust designs.


\section*{Acknowledgements}
The authors would like to acknowledge support of AFOSR Grant FA9550-15-1-0276. We also want to thank Dr. Xevi Roca and Professor Anthony Patera for useful conversations, suggestions and comments.

\bibliography{mainbib}
\bibliographystyle{elsarticle-harv}

\appendix
\section{Stochastic modeling for geometric variability}\label{ap:map}

\begin{figure}[h!]
\centering
\subcaptionbox{Reference domain (left) and physical domain (right) of square subdomain in polar coordinates. \label{fig:mapping}}
[7.5cm]{\includegraphics[scale = .9]{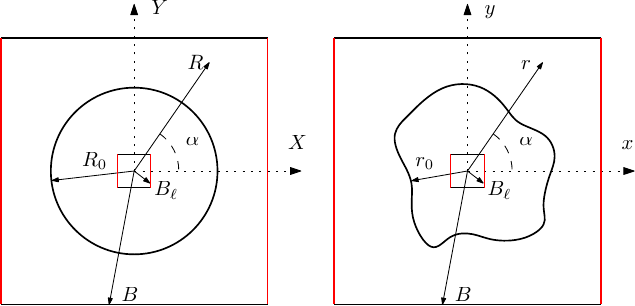}}
\hfill\subcaptionbox{Realizations of radius variation (solid) and reference radius (dashed). \label{fig:KLrad}}
[6cm]{\includegraphics[scale = 0.18]{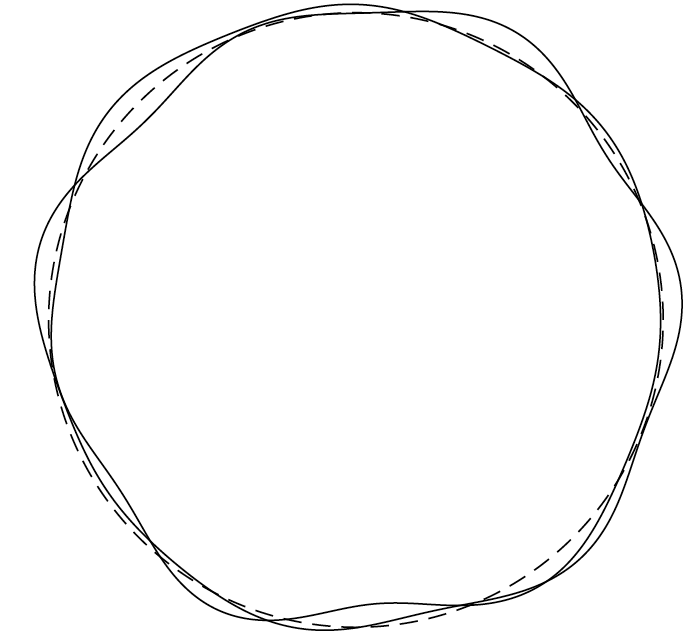}}
\caption{Mapping quantities and sample radius modifications.}
\end{figure}

In this appendix, we present the details of the deformation mapping that is employed for the circular rods/holes of the photonic crystal. We show only the case for square periodicity, as the hexagonal periodicity is a straightforward extension. For this study, we will assume that geometry variability is restricted to the radial direction and modeled with a random field of the angular coordinate $\alpha \in [0,\,2\pi]$  with known mean and covariance kernel \cite{chen2014comparison}. In particular, the geometry of the rods is characterized by the following truncated Karhunen-Lo\`{e}ve expansion of the radius:
\begin{equation}\label{eq:KL}
r_0 (\bm z)  = R_0 +\delta R_0(\bz) = R_0 \lb 1 + z_1\sqrt{\lambda_0/2} + \sum_{d=1}^{D/2} \sqrt{\lambda_d}\lb z_{2d}\sin \lp d(\alpha + \pi/2) \rp + z_{2d+1}\cos \lp d(\alpha + \pi/2)\rp  \rb \rb .
\end{equation}  
Here $\sqrt{\lambda_d} = \sigma\lp \sqrt{\pi} L_c\rp^{1/2}\exp \lp -(d\pi L_c)^2/8 \rp, d = 0,\ldots,D,$ where $\sigma$ is the variance of the covariance kernel and $L_c$ is the correlation length, which is inversely related to the decay of the KL modes. Fig. \ref{fig:KLrad} shows some realizations of  (\ref{eq:KL}). 

We now describe the mapping $\mathfrak{G}^m$ used for subdomains containing a dielectric rod. We consider square subdomains as shown in Fig. \ref{fig:KLrad}, and require that the subdomain boundaries remain fixed, as well as a small box surrounding the origin (to prevent the mapping becoming singular).  Since we are only considering deformation in the radial direction, the mapping can be expressed as $(x,y) = \mathfrak{G}^m(\refX)$ with $x  = r(\refX)\cos\alpha$ and $y = r(\refX)\sin\alpha$. 


In order to define $r(\refX)$ for the interior points in the domain, we assume that the perimeter of the rod is deformed according to \eqref{eq:KL}, and that the remaining points are linearly deformed according to 
\begin{equation}
 r(\refX) = \begin{cases} B_\ell + (r_0-B_\ell)\cdot \dfrac{R-B_\ell}{R_0 - B_\ell}, & \mbox{if } R \le R_0 \\  r_0 + (B - r_0)\cdot \dfrac{R-R_0}{B - R_0}, & \mbox{if } R > R_0  \end{cases}.
\end{equation}
The derivatives $\partial r/\partial x_i$ may be computed with the chain rule. Note there are two distinct types of boundaries, see Fig. \ref{fig:mapping}, and its distance to the origin can described as $B = \sqrt{1+(Y/X)^2}$ (red) and $B = \sqrt{1+(X/Y)^2}$ (black). The distance $B_\ell$  to the boundary of the origin box is similarly computed.

\end{document}